\documentclass[12pt,reqno]{amsart}
\usepackage{latexsym}
\usepackage{amssymb}
\usepackage[dvips]{graphicx}
\usepackage{epsfig}
\usepackage{rotating}
\usepackage{amsfonts}
\usepackage{amsmath}

\usepackage{pstricks}

\newtheorem{theorem}{Theorem}[section]

\newtheorem{lemma}[theorem]{Lemma}

\newtheorem{proposition}[theorem]{Proposition}

\theoremstyle{definition}
\newtheorem{definition}[theorem]{Definition}

\theoremstyle{remark}
\newtheorem{remark}[theorem]{Remark}

\newcommand{\Gsub}{\underline{\Gamma}}

\def\KK{{\mathbb K}}
\def\NN{{\mathbb N}}
\def\CC{{\mathbb C}}
\def\RR{{\mathbb R}}
\def\NN{{\mathbb N}}

\newcommand{\beqa}{\begin{eqnarray}}
\newcommand{\eeqa}{\end{eqnarray}}

\newcommand{\cA}{{\mathcal {A}}}
\newcommand{\cC}{{\mathcal {C}}}
\newcommand{\cB}{{\mathcal {B}}}
\newcommand{\cH}{{\mathcal {H}}}

\newcommand{\id}{{\mathrm{id}}}

\newcommand{\ba}{\begin{array}}
\newcommand{\ea}{\end{array}}

\newcommand{\veps}{\varepsilon}

\def\>{\rangle}
\def\<{\langle}

\begin{document}

\title{Some combinatorial aspects of quantum field theory}

\author{Adrian Tanasa}

\address{
LIPN, UMR CNRS 7030\\
Institut Galil\'ee \\
Univ.\ Paris 13, Sorbonne Paris Cit\'e\\
99, avenue Jean-Baptiste Cl\'ement,\\
93430 Villetaneuse, France\newline\indent
Horia Hulubei National Institute for Physics and Nuclear Engineering,\\
P.O.B. MG-6, 077125 Magurele, Romania
}

\date{}

\begin{abstract}
In this survey we present the appearance of some combinatorial notions in quantum field theory. We first focus on graph polynomials (the Tutte polynomial and its multivariate version) and their relation with the parametric representation of the commutative $\Phi^4$  field theory. We then generalize this to ribbon graphs and present the relation of the Bollob\'as--Riordan polynomial with the parametric representation of some  $\Phi^4$  field theory on the non-commutative Moyal space. 
We also review the role played by the Connes--Kreimer Hopf algebra as the combinatorial backbone of the renormalization process in field theories. We then show how this generalizes to the scalar  $\Phi^4$  field theory implemented on the non-commutative Moyal space.
Finally, some perspectives for the further generalization of these tools to quantum gravity tensor models are briefly sketched.
\end{abstract}

\maketitle

\thispagestyle{myheadings}
\font\rms=cmr8 
\font\its=cmti8 
\font\bfs=cmbx8

\markright{\its S\'eminaire Lotharingien de
Combinatoire \bfs 65 \rms (2012), Article~B65g\hfill}
\def\thepage{}

\section{Introduction\label{int} }

Within the framework of quantum field theoretical models, it has been shown 
recently that the role of various combinatorial notions is a crucial one.

Inside the area of theoretical physics, quantum field theory (QFT) represents an important part; its mathematical formalism was proven successful not only in elementary particle physics (the celebrated Standard Model) but also in condensed matter physics.

Like in other domains of theoretical physics (exact solutions of statistical-mechanical problems, random matrix models, integrable systems and so on), also in QFT  a key role is played by combinatorics. Thus, perturbative QFT relies on graph theory, 
namely on (Feynman) graphs which appear in the expansion with their combinatorial weights; moreover, as we will present in this survey, graph polynomials can be naturally related to polynomials appearing in some representation of the Feynman amplitudes in QFT. Furthermore, renormalization, which lies at the heart of QFT, can be described through some appropriate combinatorial Hopf algebras, as we will see in the sequel.

On the other hand, techniques of analytic combinatorics (such as  the Mellin transform or the saddle point method) are also used in QFT for computing and analyzing quantities like Feynman integrals (which are associated to Feynman graphs).

It is worth emphasizing that, in the last years, both combinatorialists and theoretical physicist in general (or field theorists in particular) became more and more aware of this strong relation between their domains. 
This is a natural tendency, since the unfolding of new ideas in physics is {often} tied with the development of new combinatorial methods. 
We can thus speak nowadays of an emerging domain of combinatorial physics.

\bigskip

As already announced above, in this survey we focus on the appearance
of combinatorial notions in QFT, such as graph and map polynomials and combinatorial Hopf algebras. 

Thus, in graph theory the celebrated Tutte polynomial 
(see \cite{tutte1,tutte2}) is known to characterize in a particularly elegant way a generic graph. On the other hand, the Feynman integral of such a graph can be represented, through the parametric representations, using some polynomials in a set of parameters associated to the edges of the graph. In this survey we present the relation between the Tutte polynomial (or, more exactly, its multivariate formulation \cite{sokal}) and the polynomials of this QFT parametric representation \cite{altii1, altii2, io-BR}.

\pagenumbering{arabic}
\addtocounter{page}{1}
\markboth{\SMALL ADRIAN TANASA}{\SMALL SOME COMBINATORIAL ASPECTS OF
QUANTUM FIELD THEORY}

On the other hand, we also present how the Connes--Kreimer Hopf algebra of Feynman graphs \cite{ck} encodes in a powerful manner the combinatorial structure of perturbative renormalization in QFT.

We then generalize this framework and lift it to ribbon graphs  (also known as {maps},  see for example \cite{maps} and references therein). At this level, the Tutte polynomial can be replaced in a natural manner by the Bollob\'as--Riordan polynomial \cite{br1,br2}. From the QFT point of view, the appropriate models are the non-commutative ones (we present here the scalar $\Phi^{\star\, 4}$ model implemented on the Moyal space). In this case also, a parametric representation generalizing in a highly non-trivial manner the commutative one can be defined \cite{io-BR, io-param}. One can then relate the Bollob\'as--Riordan polynomial to the polynomials of the parametric representations of these QFT models on the non-commutative Moyal space \cite{io-BR}.

Furthermore, a ribbon graph version of the Connes--Kreimer Hopf algebra of Feynman graphs can be defined \cite{io-kreimer,fab}. This algebra is proven to lie behind the renormalization of QFT models on the non-commutative Moyal space.



\section{Some algebra}
\label{sec:alg}
\renewcommand{\theequation}{\thesection.\arabic{equation}}
\setcounter{equation}{0}


In this section we briefly recall the definitions of the algebraic notions that will be used in the sequel.

\begin{definition}[\sc Algebra]\label{defn:algebra}
  A unital associative algebra $\mathcal A$ over a field $\mathbb K$ is a $\mathbb K$-linear space endowed with two algebra homomorphisms: 
  \begin{itemize}
  \item a product $m :\cA\otimes\cA\to\cA$ satisfying the \emph{associativity} condition
    \begin{equation}
      m\circ(m\otimes\id)(\Gamma)=m\circ(\id\otimes m)(\Gamma), 
\quad \text{for all
}\Gamma\in \cA^{\otimes\, 3};
\label{eq:asso}
    \end{equation}
  \item a unit $u :\KK\to\cA$ satisfying
    \begin{equation}
      m\circ(u\otimes\id)(1\otimes\Gamma)=\Gamma=m\circ(\id\otimes u)(\Gamma\otimes 1),\quad \text{for all }\Gamma\in \cA.
    \end{equation}
  \end{itemize}
\end{definition}



\begin{definition}\label{defn:coalgebra}
  A (coassociative, counital) {\it coalgebra} $\cC$ over a field $\KK$ is a $\KK$-linear space endowed with two linear homomorphisms: 
  \begin{itemize}
  \item a coproduct $\Delta :\cC\to\cC\otimes\cC$ satisfying the \emph{coassociativity} condition
    \begin{equation}
      (\Delta\otimes\id)\circ\Delta(\Gamma)=(\id\otimes \Delta)\circ\Delta(\Gamma),\quad \text{for all }\Gamma\in\cC;\label{eq:coasso}
    \end{equation}
  \item a counit $\veps :\cC\to\KK$ satisfying
    \begin{equation}
      (\veps\otimes\id)\circ\Delta(\Gamma)=\Gamma=(\id\otimes\veps)\circ\Delta(\Gamma),\quad \text{for all }\Gamma\in\cC.    \end{equation}
  \end{itemize}
\end{definition}

\begin{definition}
  A {\it bialgebra} $\cB$ over a field $\KK$ is a $\KK$-linear space endowed with both an algebra and a coalgebra structure (see Definitions~\ref{defn:algebra} and \ref{defn:coalgebra}) such that the coproduct and the counit are unital algebra homomorphisms (or, equivalently, the product and unit are coalgebra homomorphisms):
\begin{subequations}
\begin{align}
      \Delta\circ m_{\cB}&=m_{\cB\otimes\cB}\circ(\Delta\otimes\Delta),\ \Delta (1_\cB)=1_\cB\otimes 1_\cB,\\
      \veps\circ m_{\cB}&=m_{\KK}\circ(\veps\otimes\veps),\ \veps (1_\cB)=1.
      \end{align}
\end{subequations}
\end{definition}


\begin{definition}
  A {\it graded bialgebra} is a bialgebra graded as a linear space, 
  \begin{align}
      \cB=\bigoplus_{n=0}^\infty\cB^{(n)}   ,
  \end{align}
  such that the grading is compatible with the algebra and coalgebra structures:
  \begin{align}
    \cB^{(n)}\cB^{(m)} \subseteq \cB^{(n+m)}\text{ and }\Delta\cB^{(n)}\subseteq\bigoplus_{k=0}^n\cB^{(k)}\otimes\cB^{(n-k)} .
  \end{align}
\end{definition}

\begin{definition}
  A {\it connected bialgebra} is a graded bialgebra $\cB$ for which $\cB^{(0)}=u(\KK)$.
\end{definition}


\begin{definition}
  A {\it Hopf algebra} $\cH$ over a field $\KK$ is a bialgebra over $\KK$ equipped with an antipode map $S:\cH\to\cH$ obeying
  \begin{equation}
    m\circ(S\otimes\id)\circ\Delta=u\circ\veps=m\circ(\id\otimes S)\circ\Delta.
  \end{equation}
\end{definition}



A practical introduction to Hopf algebras for a combinatorial physicist can be found in \cite{ger}. 
For further details on this topic, the interested reader is referred
to \cite{Kassel} or \cite{Dascalescu}, for example.

\section{Graph theory --- the Tutte polynomial}
\label{sec:graph}
\renewcommand{\theequation}{\thesection.\arabic{equation}}
\setcounter{equation}{0}


\subsection{Some notions of graph theory}

In graph theory, one has the following definition (see, for example, \cite{mw}).

\begin{definition}
A pseudo-graph $\Gamma$ is defined as a set of vertices $V$ and a set of edges $E$ together with 
an incidence relation between them. 
\end{definition}

This means that multiple edges (i.e., edges connecting the same two vertices) and (self-)loops (i.e., edges connecting a vertex to itself) are allowed. Nevertheless, we will simply refer to these pseudo-graphs as graphs in the rest of this survey.

One can extend the definition above such
that a distinct type of edge --- the external edge --- is permitted.
Such an edge is attached to only one vertex.

In general, the terminologies 
used by graph theorists or by field theorists are different. For the sake of completeness, we present in this section both of them. Nevertheless, in the rest of this survey we use the graph theorists' language.

\begin{definition}
\noindent
\label{toate}
\begin{enumerate}
\item The number of edges at a vertex is called the {\it degree} of the respective vertex (field theorists refer to this as the {\it coordination number} of the respective vertex).
\item An edge whose removal increases the number of connected components of the respective graph is called a  {\it bridge} (field theorists refer to this as a {\it $1$-particle reducible} edge).
\item A connected subset of edges and vertices, the number of edges being the same as the number of vertices, which cannot be disconnected by removing any of the edges is called a {\it cycle} (field theorists refer to this as a {\it loop}).
\item A  {\it (self-)loop} (see above) is called a {\it tadpole} edge in QFT terminology.
\item An edge which is neither a bridge nor a self-loop is called {\it regular}.
\item An edge which is not a self-loop is called {\it semi-regular}.
\item A graph with no cycles is called a {\it forest}.
\item A connected forest is called a {\it tree}. A spanning tree is a tree connecting all the vertices of the graph.
\item 
A spanning forest with two connected components is called 
(in QFT terminology) a {\it two-tree}.
\item The {\it rank} of a subgraph $A$  is defined by
\begin{equation}
r(A):=|V|-k(A),
\end{equation}
where $k(A)$ is the number of connected components of the subgraph $A$.
\item The {\it nullity} (or {\it cyclomatic number}) of a subgraph $A$ is defined by
\begin{equation}
n(A):=\vert A\vert - r(A).
\end{equation}
\end{enumerate}
\end{definition}

\begin{figure}[bth]
\centerline{\includegraphics[width=8cm]{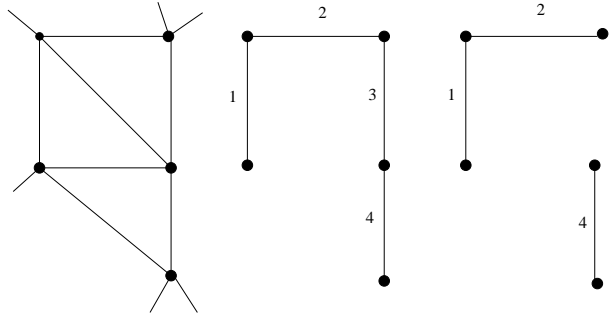}}
\caption{An example of a graph (with seven edges and six external edges). We chose a spanning tree and we label its edges by $1,2,3,4$. 
The set $\{1,2,4\}$ is a two-tree; one has two connected components (the first one formed by the edges $1$ and $2$ and the second one formed by the edge $4$). The external edges are attached to one of these two connected components.
\label{2}}
\end{figure}

We illustrate this notion by the example of Figure~\ref{2}.

\subsection{Graph polynomials and the (multivariate) Tutte polynomial}

One can define two natural operations for an arbitrary edge $e$ of some graph $\Gamma$:
\begin{itemize}
\item the {\it deletion}, which leads to a graph denoted by $\Gamma-e$,

\item the {\it contraction}, which leads to a graph denoted by $\Gamma/e$. This operation identifies the two vertices $v_1$ and $v_2$ at the ends of $e$ 
into a new vertex $v_{12}$, attributing all the edges attached to $v_1$ and $v_2$ to $v_{12}$, and then it removes $e$.
\end{itemize}

\begin{remark}
If $e$ is a self-loop, then $\Gamma/e$ is the same as $\Gamma-e$.
\end{remark}

\begin{figure}
\begin{center}
\includegraphics[scale=0.9,angle=-90]{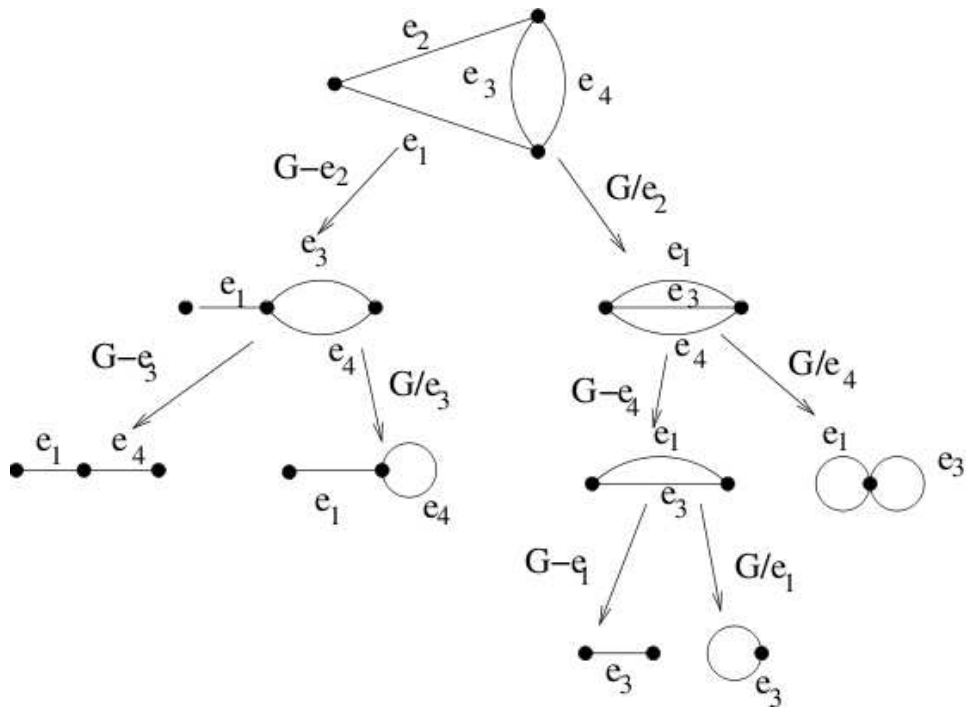}
\caption{The deletion/contraction of some graph. One is left with various possibilities (here five) of terminal forms (that is, graphs with only bridges or self-loops).}
\label{fig:cond}
\end{center}
\end{figure}

For an illustration of these two operations, see the example in Figure~\ref{fig:cond}.

\medskip

Let us now give a first definition of the Tutte polynomial.

\begin{definition}
\label{def-tutte}
 If $\Gamma$ is a graph, then its Tutte polynomial $T_\Gamma(x,y)$ 
is defined by
\begin{equation}
\label{tutte}
T_\Gamma (x,y):=\sum_{A\subset E}      (x-1)^{r(E)-r(A)} (y-1)^{n(A)}.
\end{equation}
\end{definition}

A fundamental property of the Tutte polynomial is 
the deletion/contraction relation.

\begin{theorem}
\label{defdelcontr1}
If\/ $\Gamma$ is a graph, and $e$ is a regular edge, then
\begin{equation}
T_\Gamma (x,y)=T_{\Gamma/e} (x,y)+T_{\Gamma-e} (x,y).
\end{equation}
\end{theorem}

Let us remark that this property of the Tutte polynomial is often used as its definition. However, for this property to become a mathematical definition, one needs to complete it by giving the form of the Tutte polynomial on terminal forms. Namely, for graphs with only irregular edges,  $m$ bridges and $n$ self-loops, the Tutte polynomial is given by
\begin{equation}
T_\Gamma(x,y):=x^m y^n .
\end{equation}

\medskip

A multivariate version of the Tutte polynomial exists in the literature \cite{sokal} (see also \cite{j0A,j0B}). The main idea is that, instead of having a single variable, $y$, for the number of edges, one introduces a set of variables $\beta_1,\beta_2,\ldots, \beta_{|E|}$, one for each edge. This leads to the 
following definition.

\begin{definition}
 If $\Gamma$ is a graph, then its multivariate Tutte polynomial 
is defined by
\begin{equation}
\label{multivartut}
Z_\Gamma (q,\{\beta\}):=\sum_{A\subset E}   q^{k(A)} \prod_{e\in A}\beta_e
\,.
\end{equation}
\end{definition}

Similarly, one can prove that the multivariate Tutte polynomial \eqref{multivartut} satisfies a  deletion/contraction relation, for any edge $e$. The definition of the polynomial on terminal forms (graphs with $v$ isolated vertices) is:
\begin{equation}
Z_\Gamma(q,\{\beta\}):= q^v.
\end{equation}

Through direct inspection, one can prove the 
following relation between the Tutte polynomial \eqref{tutte} and its multivariate counterpart \eqref{multivartut}:
\begin{equation}
\big[ q^{- V }   Z_\Gamma(q,\beta) \big] \Big\vert_{ \beta_e = y-1, q = (x-1)(y-1)} =
(x-1)^{k(E)  -|V|}  T_\Gamma(x,y).
\end{equation}

\bigskip

Leaving aside the Tutte polynomial, several graph polynomials have been defined and intensively studied in the literature. 
The Tutte polynomial
is a two-variable polynomial 
with
one-variable specializations such as
the chromatic polynomial or  the flow polynomial.

The {\it chromatic polynomial\/} is a graph polynomial $P_\Gamma (k)$ ($k\in \NN^\star$) which counts the number of distinct ways to color the graph $\Gamma$ with $k$ or fewer colors, where colorings are considered as distinct even if they differ only by permutation of colors. For a connected graph, this polynomial is related to the Tutte polynomial \eqref{tutte} by the relation
\begin{equation}
P_\Gamma(k)=(-1)^{|V|-1}k T_\Gamma (1-k,0).
\end{equation}

In order to define the flow polynomial, we need a finite Abelian group $G$. 
One can choose arbitrarily an orientation for each edge of the graph $\Gamma$, the result being independent of this choice. (The same type of situation appears when computing Feynman integrals, see the next section.) A {\it $G$-flow} on $\Gamma$ is a mapping 
\begin{equation}
\psi:E\to G
\end{equation}
that satisfies current conservation at each vertex. A $G$-flow on $\Gamma$ is said to be {\it nowhere-zero} if $\psi(e)\neq 0$ for all $e$. Let $F_\Gamma (G)$ be the number of nowhere-zero $G$-flows on $\Gamma$. One can prove that that this number depends only on the order $k$ of the group $G$; it can thus be written $F_\Gamma (k)$ --- it is the restriction to non-negative integers of a polynomial in $k$, the {\it flow polynomial}.
One can show that
\begin{equation}
F_\Gamma (k)= (-1)^{|E|-|V|+1}T_\Gamma (0,1-k).
\end{equation}

\section{QFT, Feynman integrals and their parametric representation}
\label{sec:qft}
\renewcommand{\theequation}{\thesection.\arabic{equation}}
\setcounter{equation}{0}

In this section we introduce a few notions of QFT, the latter being a general framework lying at the very heart of fundamental physics. QFTs give a quantum description of particles and interactions, which is naturally compatible with Einstein's theory of special relativity. As already mentioned in the Introduction, QFTs led to the Standard Model of elementary particle physics, which is one of the best experimentally tested theories (for example in huge particle accelerators like the CERN's LEP). QFT's mathematical formalism also successfully applies to other branches of theoretical physics, like condensed matter or statistical physics.

The interested reader may consult any of the very good textbooks on QFTs, such as 
\cite{iz} or \cite{peskin}.

\bigskip

The simplest field theoretical model is the {\it $\Phi^4$ scalar model}; it consists of a single type of field $\Phi(x)$, a scalar field,
\begin{equation}
\Phi\; : \; \RR^4 \; \to\; \KK,
\end{equation}
where the target space $\KK$ is either $\RR$ or $\CC$. 

Note that we use here $\RR^4$, which corresponds to a Euclidean metric. If one needs to work with a Minkowskian metric, then $\RR^4$ needs to be replaced by $\RR^{1,3}$ (the time being singled out from the three spatial dimensions). The number $4=1+3$ corresponds to the fact that we work on a four-dimensional space-time.

For a field theoretical model to be defined, one needs to write down an {\it action}, which is, from a mathematical point of view, a functional of the fields of the model. For the real $\Phi^4$ model, the action reads
\begin{equation}
\label{act}
S[\Phi (x)]=\int_{\RR^4} d^4 x \left[ 
\frac{1}{2} \sum_{\mu=1}^4\partial_\mu \Phi(x) \partial_\mu
\Phi(x) +\frac{1}{2}m^2{\Phi(x)}^2 + V_{\text{int}}[\Phi(x)]
\right],
\end{equation} 
where the real parameter $m$ represents the mass. Furthermore, we have $\partial_\mu=\frac{\partial}{\partial x_\mu}$, and $V_{\text{int}}[\Phi(x)]$ stands for the {\it interacting potential}, namely $  \frac{\lambda}{4!} {\Phi (x)}^4$ ($\lambda$ being a real parameter, called the coupling constant).

The formula \eqref{act} is written in {\it configuration space} (or {\it direct space}). By applying the Fourier transform, one can write it in {\it momentum space} (which is actually the space where the elementary particle physics computations are usually performed):
\begin{equation}
\label{act-p}
\tilde S[\tilde \Phi]=\int_{\RR^4} d^4 p \left[ 
\frac{1}{2} \sum_{\mu=1}^4 p_\mu  p_\mu \tilde \Phi^2
+\frac{1}{2}m\tilde \Phi^2 + \tilde V_{\text{int}}[\tilde \Phi]
\right].
\end{equation}

Let us now further analyze this action. One has a quadratic part in $\Phi$ and a non-quadratic one. The quadratic part corresponds to the {\it propagation} (the {\it free theory}), while the non-quadratic part corresponds to the {\it interaction}. 




The  functional integration is introduced as the product of integrals at each space point $x$ (multiplied by some irrelevant normalization factor), 
${\mathcal D}\phi (x):={\mathcal N}\prod_x \int d\Phi(x)$.
Nevertheless, such an infinite product of Lebesgue measures is mathematically ill-defined. For a well-defined QFT measure, the interested reader may consult, for example, the review article \cite{review-riv}.

The definition of the {\it partition function} is
\begin{equation}
\label{defz}
Z^{\mathrm{}}:=\int {\mathcal D}\Phi (x) e^{-S[\Phi(x)]}.
\end{equation}
The physical information of a theory is encoded in  {\it $n$-point functions} (or 
{\it correlation functions}) which are defined by
\begin{align}
\nonumber
S^{(N)}(x_1,x_2,\ldots,x_N)&:=\frac{1}{Z}\int  {\mathcal D}\phi (x)
\Phi (x_1)\Phi(x_2)\cdots \Phi(x_N)e^{-S[\Phi]}\\
&=\langle \Phi(x_1)\Phi(x_2)\ldots\Phi(x_N)\rangle.
\label{defS}
\end{align}

In general,
 one is unable to find exact expressions for these correlation functions. 
In such cases, the tool used in theoretical physics is the {\it perturbative expansion}, i.e., an {\it expansion of the exponential} above in powers of $\lambda$. The coefficients of such an expansion are sums of multiple integrals (see \eqref{act} and \eqref{defS});
the number of these integrals grows rapidly with increasing {\it order} in perturbation theory (i.e., power of the coupling constant $\lambda$).

\begin{figure}
\begin{center}
\includegraphics[scale=0.7,angle=-90]{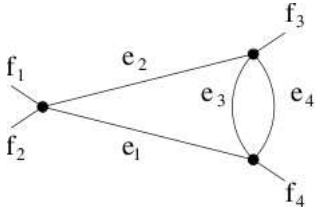}
\caption{ A $\Phi^4$ Feynman graph, with four internal edges ($e_1,e_2,e_3, e_4$) and four external edges  ($f_1,f_2,f_3, f_4$). It is a graph at the third order in perturbation theory (it has three vertices).}
\label{fig:nor}
\end{center}
\end{figure}

To these multiple integrals, called {\it Feynman integrals}, 
one associates {\it Feynman graphs}, which are very useful in the organization of the expansion coefficients. For the $\Phi^4$ model exhibited here, the graphs have degree four at each vertex; moreover one has internal and external edges (see, for example, Figure~\ref{fig:nor}).

Each of these graphs comes with its combinatorial weight (which is highly non-trivial because of the non-labeling of the edges). 
Furthermore, the Feynman integrals can then be manipulated using various 
techniques of analytic combinatorics (the Mellin transform or the saddle-point approximation).

In momentum space, one orients each edge $e$ (internal or external) and associates to it  some momentum $p_e$. 
The contribution from any internal edge comes from the quadratic part of the action \eqref{act-p} and reads
\begin{equation}
\label{propa}
C(p_e)=\frac{1}{p_e^2+m^2}.
\end{equation}
This contribution is called the {\it propagator}.

The contribution of some vertex is given by the coupling constant $\lambda$ times a $\delta$-function of conservation of the incoming/outgoing momenta at the respective vertex. 
If the momenta of the external edges are fixed (they correspond to physical states), 
the momenta of the internal edges are free and one has to integrate over
them, leading to
the respective Feynman integral.

For the example of Figure~\ref{fig:nor}, all this reads
\begin{multline}
\label{ex}
\lambda^3\int \left(\prod_{i=1}^4 d^4 p_{e_i}\right)
\left(\prod_{i=1}^4 \frac{1}{p_{e_i}^2+m^2}\right)\\
\cdot
\delta(p_{f_1}+p_{f_2}-p_{e_1}-p_{e_2})
\delta(p_{e_1}+p_{e_3}-p_{e_4}-p_{f_4})
\delta(p_{e_2}-p_{e_3}+p_{e_4}-p_{f_3}).
\end{multline}
We leave it as an exercise to the interested reader to find out the orientation of the edges which was chosen for such a Feynman integral to occur.

Because of the presence of the three $\delta$-functions in \eqref{ex}, the number of remaining integrals is equal to two, which is actually the number of independent cycles of the Feynman graph (this being a general result in QFT).

At the end of the story, for the Feynman integral \eqref{ex} all this leads to a logarithmic  divergence in the high energy ($|p|\to \infty$) regime of the internal momenta (the so-called {\it ultraviolet regime}).
The appearance of such divergences is in fact a very frequent phenomenon in QFT; it is the {\it renormalization} process which, when possible, takes care of these infinities in a highly non-trivial way (see, for example, the book \cite{book-riv}).

\bigskip

Let us now give some insight in this process of {renormalization} in QFT. 
If the considered QFT model is renormalizable, 
the graphs which lead to divergences should correspond to terms present in the action. 
In order to illustrate what we mean by this, let us go back to the example of the $\Phi^4$ model. The graphs which lead to the various divergences of the model need to have two or four external edges. Thus, these divergences can be ``cured'' by an appropriate renormalization of the parameters of the action \eqref{act} 
(for example, the mass $m$ and the coupling constant $\lambda$). 
For instance, a graph with four external edges (each of these four external edges being
 associated to a field $\Phi$) corresponds 
to the renormalization of the coupling constant $\lambda$. 
This comes from the fact that the coupling constant $\lambda$ is indeed
the parameter by which the  $\Phi^4$ term in the action  \eqref{act}
is multiplied.

\begin{figure}[bth]
\centerline{\includegraphics[width=4cm]{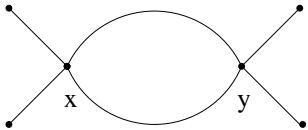}}
\caption{A four-point Feynman graph. It has two vertices, localized, in configuration space, in $x$ and $y$. 
\label{locality-xy}}
\end{figure}

Let us take as an example the graph $\Gamma$ with four external edges in Figure~\ref{locality-xy}. 
In configuration space, the Feynman integral to renormalize reads
(up to irrelevant normalization factors)
\begin{equation}
\label{ea}
\phi (\Gamma)=\int dx dy (C(x,y))^2 C(x,\bar y_1)C(x,\bar y_2)C(y,\bar y_3)C(y,\bar y_4),
\end{equation}
where we have denoted the four external points that the external edges hook to by $\bar y_1,\bar y_2,\bar y_3 , \bar y_4$.
Note that the propagators $C(x,y)$ are given by the inverse Fourier transforms of the momentum space propagators \eqref{propa} --- they are associated to the quadratic part of the configuration space action \eqref{act}, i.e.,  they are integral representations of heat kernel functions (by a slight abuse of notation we have denoted propagators in both configuration and momentum space by $C$).

We now subtract a term corresponding to the configuration space region 
\begin{equation}
\label{loc}
x\sim y. 
\end{equation}
The integral \eqref{ea} is divergent and can be rewritten, using an 
appropriate Taylor expansion, as
\begin{multline}
\label{taylor}
\int dx dy (C(x,y))^2 C(x,\bar y_1)C(x,\bar y_2)\\
\cdot
\left( C(x,\bar y_3)C(x,\bar y_4)
+ \int_0^1 dt (y-x) \nabla \left( C(x+t(y-x),\bar y_3)C(x+t(y-x),\bar y_4) \right)\right).
\end{multline}
One can prove that the second term in the sum above --- the renormalized Feynman integral --- is finite (see, for example, the book \cite{book-riv} for details), while the first is the one containing the divergence --- it will be subtracted. This term is called {\it local counterterm}:
\begin{equation}
\label{ct}
\tau_\Gamma  \phi (\Gamma)=
\int dx(\Phi (x))^4\int dy (C(x,y))^2.
\end{equation}
The adjective ``local" comes from the fact that 
this expression 
is associated, as shown above, to the configuration space region \eqref{loc}. 
The operator $\tau$ 
represents the first term in the Taylor expansion \eqref{taylor}.
Note also that in \eqref{ct} we have explicitly considered 
the multiplication with the four fields $\Phi(x)$   --- 
this divergent term can thus be reabsorbed, as already mentioned above,
 in a {redefinition of the coupling constant $\lambda$}.
Since the model \eqref{act} is translation-invariant, the second integral in \eqref{ct} is independent of $x$, and thus the counterterm \eqref{ct} is of the same form as the interaction term in the action.
We have thus shown that {\it the subtraction of a local counterterm makes the renormalized Feynman integral finite.}

The locality of the counterterms can be seen as 
{\it sending all of the external edges to the same point\/} (here $x$), 
again, in configuration space.  
The phenomenon described here can also be exhibited in momentum space 
(see, for example, the review paper \cite{review-riv} for the same example treated in momentum space). 

Nevertheless, subtraction of the divergences can be more involved. 
One can have graphs with subdivergences.  
In order to introduce the general (Bogoliubov) subtraction operator for a Feynman graph, 
we 
first need the following definition.

\begin{definition}
 A {\it Zimmermann forest} $\mathcal F$ of a Feynman graph $\Gamma$ is a set of subgraphs of $\Gamma$ such that
\begin{equation}
 \gamma\cap \gamma'=\emptyset \text{ or } (\gamma\subset\gamma' \text{ or } \gamma'\subset\gamma), \quad \text{for all } \gamma,\gamma'\in {\mathcal F}. 
\end{equation}
\end{definition}

\begin{remark}
 The Zimmermann forest defined above is a notion distinct from that of a forest in graph theory (see Definition~\ref{toate}(7)). 
\end{remark}

The (Bogoliubov) subtraction operator reads
\begin{equation}
\label{R}
\bar R:=\sum_{{\mathcal F}}\prod_{\gamma \in {\mathcal F}} (-\tau_ \gamma),
\end{equation} 
where the sum is over all Zimmermann forests $\mathcal F$ of superficially divergent subgraphs 
$\gamma$ (including the empty forest). The superficially divergent graphs are, in the case of the 
$\Phi^4$ model, the graphs with two or four external edges (this result comes from the so-called 
power counting theorem).

Applying this subtraction operator to a Feynman integral $\phi(\Gamma)$ extracts the divergences. The remaining part, the {\it renormalized Feynman integral} 
\begin{equation}
\phi_{+}(\Gamma)=\bar R (\phi (\Gamma)),
\end{equation}
 is finite.

\begin{figure}
\centerline{\epsfig{figure=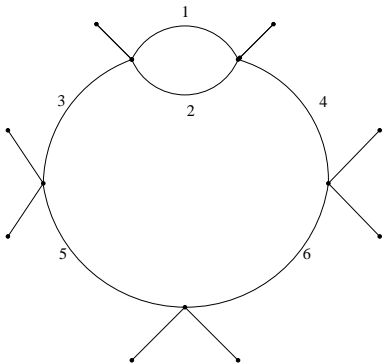,width=5cm} }
\caption{A Feynman graph, with six internal edges and with a subdivergence given by the subgraph made of the internal edges $1$ and $2$.}
\label{big}
\end{figure}

An example of a Feynman graph with subdivergences is given in Figure~\ref{big}, where the Feynman graph of Figure~\ref{locality-xy} is a subgraph. 
Renormalizing this subdivergence leads to the renormalization of the whole graph. This can be explicitly seen from the following facts.
The  Zimmermann forests of  this graph are
\begin{equation} 
\label{listz}
\emptyset, \{ \gamma \}. 
\end{equation}
The renormalized Feynman integral thus reads
\begin{equation}
\label{ex-qft}
\phi_{+}(\Gamma)=\bar R (\phi (\Gamma))
=\phi(\Gamma) - (\tau_\gamma(\phi(\gamma)))\phi(\Gamma/\gamma),
\end{equation}
where we have used the general property $\tau_\gamma \phi (\Gamma) = \phi (\Gamma/\gamma)\tau \phi (\gamma)$, $\gamma$ being a subdivergence of the Feynman graph $\Gamma$.
Note that we wrote $\gamma/\gamma'$ for the {\it cograph} obtained by shrinking the subgraph $\gamma'$ inside $\gamma$ --- {\it shrinking} a subgraph means to replace its internal structure by a point.

This graph-by-graph renormalization method is known under the name of the 
Bo\-go\-liu\-bov--Parasiuk--Hepp--Zimmermann
 (BPHZ) renormalization scheme, and it is one of the most frequently used in QFT.


For a general Feynman graph, it is also 
possible to have {\it overlapping divergences}, 
i.e.,  divergent subgraphs with a non-empty intersection. 
In this case, one can use 
a classification theorem of Zimmermann forests in order to prove finiteness 
of renormalized Feynman amplitudes.
In this case, one needs 
a scale decomposition of the graph's propagators. 
This means that each propagator $C$ is written as a sum of distinct contributions,
each such contribution corresponding to a distinct energy scale
(see the book \cite{book-riv} for a detailed presentation of these issues).

\medskip

Let us also mention that, once the renormalization techniques are performed, the renormalized (and hence finite) integral leads to physical quantities which are measured in elementary particle collider experiments with an extremely high accuracy.

\bigskip

In the case presented here, the $\Phi^4$ model, the mass $m$ prevents the integral to be divergent also in the {\it infrared regime} (that is, for $|p|\to 0$; see, for example equation~\eqref{ex}, which is convergent in this regime).

\medskip

 For the sake of completeness, let us also mention that, when computing Feynman integrals, 
mathematical physicists are interested in graphs with no bridges 
(called $1$-particle irreducible ($1$PI) graphs). 
This comes from the fact that the momentum of the bridge cannot flow towards 
its high energy  regime, because it is directly related to the external (and thus always finite) momenta. These graphs are thus of no interest for a renormalization analysis.

\medskip

Let us now proceed further and introduce the {\it parametric representation} of a Feynman integral. The main idea is to write each of the internal propagators \eqref{propa} of the integral as an integral over some parameter $\alpha$:
\begin{equation}
\label{trick}
C(p_e)=\int_0^\infty d\alpha_e e^{-\alpha_e (p_e^2+m^2)},\quad  e=1,2,\ldots , E.
\end{equation}
Inserting these formulas in the general expression of a Feynman amplitude allows to 
integrate out (through Gaussian integrations)
the internal momenta $p_e$. 

The Feynman integral then becomes
\begin{equation} 
\phi(\Gamma)=\int_0^{\infty} 
\frac{e^{- V_\Gamma(p_{{\rm ext}},\alpha)/U_\Gamma (\alpha) }}{U_\Gamma (\alpha)^{2}} 
\prod_{\ell=1}^{|E|}
( e^{-m^2 \alpha_\ell} 
d\alpha_\ell 
),
\end{equation}
where $U_\Gamma(\alpha)$ is a polynomial in the set of parameters $\alpha$ and 
$V_\Gamma(p_{{\rm ext}},\alpha)$ is a polynomial in the set of external momenta $p_{{\rm ext}}$ and  the set of parameters $\alpha$. 
It can be shown that these polynomials depend only on the underlying graph.
More specifically, one has
\begin{equation}
\label{s1}
U_\Gamma (\alpha)= \sum_{\mathcal T} \prod_{\ell \not \in {\mathcal T}} \alpha_\ell
\end{equation}
and 
\begin{equation}
\label{s2}
V_\Gamma (p_{{\rm ext}},\alpha)= \sum_{{\mathcal T}_2} \prod_{\ell \not \in {\mathcal T}_2} \alpha_\ell  \bigg(\sum_{i \in
  E({\mathcal T}_2)} p_i\bigg)^2, 
\end{equation}
where the first sum is over all trees $\mathcal T$ 
of the graph $\Gamma$, while the second sum is over two-trees
${\mathcal T}_2$,
which, as already stated in the previous section 
(see Definition~\ref{toate}(9)) 
separates the graph 
in two connected components and 
one can identify the external edges which connect to one or the other of these connected components.
The symbol $E({\mathcal T}_2)$ denotes 
one of the connected components thus obtained.

\begin{remark}
By momentum conservation, the total  momentum of one of these connected components (for example  $E({\mathcal T}_2)$) is equal to the total momentum of the other connected component. 
\end{remark}

For the example of Figure~\ref{fig:nor}, one has
\begin{align}
\nonumber
U_\Gamma(\alpha)&=\alpha_3\alpha_4+\alpha_2\alpha_4+\alpha_2\alpha_3+\alpha_1\alpha_3+\alpha_1\alpha_4,\\
V_\Gamma(p_{{\rm ext}},\alpha)&=(p_{f_1}+p_{f_2})^2\alpha_1\alpha_2(\alpha_3+\alpha_4)+p_{f_4}^2\alpha_1\alpha_3\alpha_4+p_{f_3}^2 \alpha_2\alpha_3\alpha_4.\nonumber
\end{align}


\section{Relation between the multivariate Tutte polynomial and the polynomials of the parametric representation}
\label{sec:rel}
\renewcommand{\theequation}{\thesection.\arabic{equation}}
\setcounter{equation}{0}

Let us now exhibit the relation between the notions presented in the two previous sections, the multivariate Tutte polynomial on the one hand, 
and the polynomials of the parametric representation of Feynman integrals
on the other hand.

We start with a deletion/contraction relation for the polynomials
$U_\Gamma(\alpha)$.

\begin{theorem}
\label{delconU}
For any semi-regular edge $e$, we have
\begin{equation} 
\label{delcontr1}
 {U}_\Gamma (\alpha) = \alpha_e \;  {U}_{\Gamma-e} (\alpha) + {U}_{\Gamma/e} (\alpha).
\end{equation}
Moreover, the terminal form evaluation is
\begin{equation}\label{firsttermin} 
{U}_\Gamma (\alpha)  = \prod_{e}  \alpha_e,
\end{equation}
for $G$ consisting only of self-loops attached to isolated vertices.
\end{theorem}

In order to prove this result, we shall 
use the Grassmann representation 
(see Appendix~\ref{grass}) of the polynomial ${U}_\Gamma $. 

A key ingredient for this representation is the {\it incidence matrix} 
$(\epsilon_{e,v})$ of the graph. In this matrix, the rows are indexed
by the edges $e$ of the graph, $e=1,2,\ldots E$, and the columns are indexed
by the vertices $v$ of the graph, $v=1,2,\ldots,V$.
By definition, the entry $\epsilon_{e,v}$ is equal to $1$ if the edge $e$
is outgoing from the vertex $v$, 
$-1$ if the edge $e$ is incoming at the vertex $v$, 
and it is $0$ if the edge $e$ is not incident to the vertex $v$. 
Let us point out that this matrix encodes the 
complete information of a Feynman graph without self-loops. 
In order to encode the information from self-loops also (if they are present), 
one needs to define an additional matrix 
$(\eta_{e,v})$, 
with non-vanishing entries $\eta_{e,v}$ if $e$ is a self-loop incident to the vertex $v$ 
 and  $\eta_{e,v}=0$ otherwise (note that $\epsilon_{e,v}=0$ if $e$ is 
a self-loop incident to the vertex $v$, see 
 \cite{io-BR} for details).

The proof of Theorem~\ref{delconU} is then based on the following 
auxiliary result.

\begin{lemma}
\label{lemma1}
The polynomial ${U}_\Gamma (\alpha)$ is given by
\begin{equation}
{\mathrm{det}}\, Q_\Gamma,
\end{equation}
where the graph matrix $Q_\Gamma$ is the
square matrix 
\begin{equation}
Q_\Gamma  =
\begin{pmatrix} \alpha_e &   (-\epsilon_{e,v} ) \\
 {}^t(-\epsilon_{e,v} ) & 0  \\ \end{pmatrix}.
\end{equation}
{\em(}$Here, {}^t(\,\cdot\,)$ 
denotes the transpose of a matrix{\em)}.
\end{lemma}
\begin{proof}
One can prove this by performing the Gaussian integrations obtained from inserting the propagator form \eqref{trick} in the general form of a Feynman integral (see again \cite{io-BR} for details).   
\end{proof}

\begin{proof}[Proof of Theorem~\ref{delconU}]
As a direct consequence of Lemmas~\ref{lemma1} and \ref{lemma2}, we have
\begin{equation}
\label{ince}
{U}_G (\alpha_e) = \int \prod_{v,e}  d\chi_v d\omega_v  d\chi_e d\omega_e
e^{ - \alpha_e\chi_e \omega_e}
e^{- \chi_e
\epsilon_{e,v} \chi_v  + \omega_e
\epsilon_{e,v} \omega_v } .
\end{equation}
Note that we have split each of the two sets of Grassmann variables ($\chi$ and $\omega$) 
required for the representation of such a determinant (see Lemma~\ref{lemma2})
in two subsets, one of cardinality $E$ (indexed by $e$) and one of cardinality $V$ (indexed by $v$); this corresponds to the block form of the complete matrix $Q_\Gamma$.

Let now $e$ be a semi-regular edge connecting vertices $v_1$ and 
$v_2$.
We have (see the definition~\eqref{defg} for a single Grassmann variable)
\begin{equation} 
\label{sumita}
e^{- \alpha_e\chi_e \omega_e} = 1 + \alpha_e\omega_e\chi_e . 
\end{equation}
We denote the two terms of this sum by $\det Q_{\Gamma,e,1}$ and by  $\det Q_{\Gamma,e,2}$, respectively. 
Let us first closely investigate the term  $\det Q_{\Gamma,e,1}$.
By direct inspection of the various terms coming from similar expansions of the other Grassmann exponentials in \eqref{ince}, one can prove (using the Grassmann integration rules \eqref{intg}) that (see again \cite{io-BR} for details)
\begin{multline}
\det Q_{\Gamma,e,1}= \int  \prod_{e'\ne e,v}  d \chi_{e'}d\omega_{e'}  d\chi_v d\omega_v   
(\chi_{v_1} - \chi_{v_2})  (\omega_{v_1} - \omega_{v_2}) 
\\
\cdot
e^{-\sum_{e'\ne  e }\alpha_e'\chi_{e'}
\omega_{e'}} 
e^{-\sum_{e' \ne e ,v}   \chi_{e'}
\epsilon_{e',v}\chi_v+ \sum_{e' \ne e ,v}   \omega_{e'}
\epsilon_{e',v}\omega_v   } .
\end{multline}
Let us now perform the following (triangular)  change of variables  (with unit Jacobian)
\begin{align}
\nonumber
\hat \chi_{v_1} &= \chi_{v_1} - \chi_{v_2}, \\  
\hat \chi_{v} &= \chi_{v},\quad \text{for all } v \ne v_1,\nonumber
\end{align}
and the same for the set of $\omega$ Grassmann variables.
This leads to the identification (see again \cite{io-BR} for details)
\begin{equation}
\det Q_{\Gamma,e,1} =  \det Q_{\Gamma/e}. 
\end{equation}
The second term of the expression \eqref{sumita} similarly leads to
\begin{equation}
\det Q_{\Gamma,e,2} = \alpha_e \det Q_{\Gamma-e}.
\end{equation}
We refer the interested reader again to \cite{io-BR} for obtaining the terminal form evaluation.
\end{proof}

\begin{remark}
The Grassmann calculus method presented here for the proof of Theorem~\ref{delconU} is not the simplest method to obtain this result (which can be obtained for example using the explicit form \eqref{s1} of the polynomial $U_\Gamma (\alpha)$).  Nevertheless, the Grassmann calculus method is particularly useful to get, for example, formula~\eqref{s1} (or formula~\eqref{s2} for the polynomial $V_\Gamma (p_{\text{ext}},\alpha)$) directly from the respective Feynman integral. Moreover, such a method can be used to obtain the parametric representation of more complicated QFT models, like for example the non-commutative scalar $\Phi^{\star\, 4}$ model (see Subsection~7.2). 
\end{remark}

\bigskip

The situation is similar for the second polynomial $V(p_{{\rm ext}},\alpha)$ (see \cite{io-BR} for details). 
Let us now show how the polynomial $ {U}_\Gamma (\alpha)$ can be obtained as a limit of the multivariate Tutte polynomial $Z_G (q,\beta)$.  
One needs to consider
\begin{equation}
q^{- k(\Gamma) }Z_\Gamma (q,\beta).
\end{equation}
Taking the limit $q \to 0$
one obtains a sum over maximally spanning subgraphs $A$,
that is, subgraphs with $k(A)=k(\Gamma)$:
\begin{equation} S_{\Gamma} (\beta)=\lim_{q\rightarrow 0}q^{- k(G) }Z_\Gamma (q,\beta)=\sum_{A
\mathrm{ \; \; maximally  \; \; spanning  \; \; }  E  } \quad
\prod_{e\in A} \beta_e .
\end{equation}

If one now retains only the terms of lowest degree of homogeneity
in $\beta$, one obtains a sum over maximally spanning graphs
with lowest number of edges, i.e., maximally spanning acyclic graphs (or 
{spanning forests}, see Definition~\ref{toate}(7))  of $\Gamma$. We denote this lowest number of edges by $p$ (note that $p=|V(\Gamma)|-k(\Gamma)$). We have
\begin{equation} F_{\Gamma} (\beta)=\sum_{{\mathcal F}
\mathrm{ \; \; maximally  \; \; spanning  \; \;  forest    \; \;   of  \; \; }  G  } \quad
\prod_{e\in {\mathcal F}} \beta_e .
\end{equation}

Finally, if one divides
$F_{\Gamma} (\beta)$ by $\prod_{e \in E} \beta_e$
and makes the change of variables $\alpha_e = \beta_e^{-1}$, 
the polynomial  $ {U}_\Gamma (\alpha)$ is obtained.

All this is summarized by the formula
\begin{equation}
U_\Gamma(\alpha)=\left[\left(\prod_{e \in E} \beta_e\right)\lim_{q'\rightarrow 0}
\frac{1}{(q')^{p(\Gamma)}}\lim_{q\rightarrow 0}q^{- k(\Gamma)}Z_\Gamma (q,q'\beta)\right]{}_{\beta_e^{-1}=\alpha_e}.
\end{equation}
Note that the variable $q'$ was introduced to take the limit retaining only the lowest degree of homogeneity in $\beta$ from $S_{\Gamma} (q'\beta)$ (see above).

\medskip

Let us also mention here the universality of the Tutte polynomial. This, also known as the ``recipe theorem'', means that any Tutte--Grothendieck invariant must be an evaluation of the Tutte polynomial, with the necessary substitutions given by the recipe (see, for example, the survey papers \cite{j0A,j0B}).
The relation between the Symanzik polynomial $U_\Gamma$ and the multivariate Tutte polynomial $Z_\Gamma$ can thus be deduced from this result and Theorem~\ref{delconU}.

\section{Combinatorial Connes--Kreimer Hopf algebra of Feynman graphs}
\label{sec:ck}
\renewcommand{\theequation}{\thesection.\arabic{equation}}
\setcounter{equation}{0}

Consider now the unital associative algebra  $\cH$ freely generated by $1$PI Feynman graphs of the $\Phi^4$ model (including the empty set, which we denote by $1$). The product $m$ is bilinear, commutative, and given by the operation of disjoint union.






Let us give a formal definition of the Feynman rules, already introduced in the previous section.

\begin{definition}
\label{def:FeynmanRules}
  The 
{\it Feynman rules} are a homomorphism $\mathbf{\phi}$ from $\cH$ to some target space $\cA$. 
\end{definition}

For the sake of completeness, let us also remark that this target algebra is naturally equipped with some Rota--Baxter algebraic structure. Recently,
the algebraic implications one has when relaxing this particular condition
were studied in \cite{patras}.

In a generic QFT model, a special role in the process of renormalization is played by the {primitively divergent graphs}. 

\begin{definition}
A {\it primitively divergent graph} of a QFT model is a graph whose Feynman integral is divergent but which does 
not contain any subgraph for which the Feynman integral is also divergent.
\end{definition}

In the particular case of the commutative $\Phi^4$ model, this class of graphs is formed by the graphs with two or four external edges that do not contain any graph of two or four external edges as subgraph. Let us also remark that the  parametric representation presented in Section~\ref{sec:qft} is one of the elegant ways to prove this type of result.

\begin{definition}
\label{def:projectionT}
  The {\it projection} $\mathbf{T}$ is a map from $\cA$ to $\cA$ which 
satisfies the following property: 
For all $\Gamma\in\cH$, $\Gamma$ primitively divergent,
  \begin{align}
    (\id_{\cA}-T)\circ\phi(\Gamma)<\infty.\label{eq:conditionT}
  \end{align}
  This means that if $\phi(\Gamma)$ is divergent then its overall divergence is completely included in $T\circ\phi(\Gamma)$.
\end{definition}

From now on we denote the set of two- and four-external edges subgraphs of a generic graph $\Gamma$ by $\Gsub$. This is also referred to as the set of superficially divergent subgraphs of $\Gamma$.

\begin{lemma}[{\sc\cite[Lemma 3.2]{fab}}]
\label{lem:coassociative1}
  Let\/ $\Gamma\in\cH$. Provided
 \begin{enumerate}
 \item for all $\gamma\in\Gsub$ and $\gamma'\in\underline{\gamma}$, 
the graph $\gamma/\gamma'$ is superficially divergent,
\label{item:condition1}
 \item for all $\gamma_{1}\in\cH$ and $\gamma_{2}\in\cH$ such that $\gamma_1$ and $\gamma_2$ are superficially divergent, there exists gluing data $G$ such that $(\gamma_{1}\circ_{G}\gamma_{2})$ is superficially divergent,\label{item:condition2}
\end{enumerate}
  the following coproduct is coassociative:
\begin{subequations}
  \label{eq:coproduct1}
  \begin{align}
\Delta\Gamma&=\Gamma\otimes 1
+1\otimes\Gamma+\Delta'\Gamma,\label{eq:Deltaprime}\\
\Delta'\Gamma&=\sum_{\gamma\in\Gsub}\gamma\otimes\Gamma/\gamma.
\end{align}
\end{subequations}
\end{lemma}

By {\it gluing data} we understand a bijection between the external edges of the graph to be inserted and the edges of the propagator (or vertex) where the insertion is done.

\medskip

\medskip

Let us now illustrate how this result fits the commutative $\Phi^{4}$ model. 
The first condition of the lemma is trivial, since when one shrinks a subgraph, the number of external edges is conserved. Moreover, since $\gamma'$ has two or four external edges, the cograph is also a $\Phi^4$ graph.
Let us now check condition \ref{item:condition2} of Lemma~\ref{lem:coassociative1}. We consider two graphs $\gamma_{1}$ and $\gamma_{2}$ with two or four external edges. We consider $\gamma_{0}=\gamma_{1}\circ_{G}\gamma_{2}$ for any gluing data $G$. Let $V_{i},\,I_{i}$ and $E_{i}$ be the numbers of vertices, internal, and external edges of $\gamma_{i},\,i\in\{0,1,2\}$, respectively. For all $i\in\{0,1,2\}$, we have
\begin{subequations}
  \label{eq:ConservationExtLines}
  \begin{align}
    4V_{i}&=2I_{i}+E_{i},\\
    V_{0}&=
    \begin{cases}
      V_{1}+V_{2},&\text{if $E_{2}=2$,}\\
      V_{1}+V_{2}-1,&\text{if $E_{2}=4$,}
    \end{cases}\\
    I_{0}&=
    \begin{cases}
      I_{1}+I_{2}+1,&\text{if $E_{2}=2$,}\\
      I_{1}+I_{2},&\text{if $E_{2}=4$,}
    \end{cases}
  \end{align}
\end{subequations}
which proves that $E=E_{1}$. Then, as soon as $\gamma_{1}$ is primitively divergent, so is $\gamma_{0}$. Concerning condition \ref{item:condition1}, note that $\gamma''=\gamma/\gamma'$ if and only if there exists $G$ such that $\gamma=\gamma''\circ_{G}\gamma'$, which allows to prove that condition \ref{item:condition1} also holds and that the coproduct (\ref{eq:coproduct1}) is coassociative.

\bigskip

 Let the coproduct $\Delta:\cH\to\cH\otimes\cH$ be defined by
\begin{equation}
\Delta \Gamma = \Gamma \otimes 1 + 1 \otimes \Gamma + \sum_{\gamma\in \Gsub} \gamma \otimes \Gamma/\gamma, \ \forall\Gamma\in\cH.\label{eq:NCcoproduct}
\end{equation}
Furthermore, let us define the counit $\varepsilon:\cH\to\mathbb K$ by
\begin{equation}
\varepsilon (1) =1,\ \varepsilon (\Gamma)=0,\quad \text{for all }\Gamma\ne 1.
\end{equation}
This means that, for any non-trivial element of $\cH$, the counit returns a trivial answer and, equivalently, only for the trivial element of $\cH$ (the empty graph $1$), the result returned by the counit is non-trivial. 
Finally the antipode is given recursively by
\begin{align}
\nonumber
  S:\cH&\to\cH\\
\nonumber
S(1)&=1,\\
  \Gamma&\mapsto-\Gamma-\sum_{\gamma\in\Gsub}S(\gamma)\Gamma/\gamma.
\label{eq:Antipode}
\end{align}

\medskip

Then the following result holds.

\begin{theorem}[{\sc\cite[Theorem~1]{ck}}]
The quadruple $(\cH,\Delta,\veps,S)$ is a Hopf algebra.
\end{theorem}

Let us also notice that $\cH$ is graded by the loop number.

\bigskip

Let now $f,g\in\text{Hom}(\cH,\cA)$, where $\cA$ is the range algebra of the projection $T$ (see above). The convolution product $\ast$ in $\text{Hom}(\cH,\cA)$ is defined by
\begin{equation}
  f\ast g=m_{\cA}\circ(f\otimes g)\circ\Delta_{\cH}.\label{eq:convolution}
\end{equation}
Let $\phi$ be the Feynman rules and $\phi_{-}\in\text{Hom}(\cH,\cA)$ the twisted antipode, defined
 recursively by 
\begin{equation}
  \phi_{-}(\Gamma)=-T\big(\phi(\Gamma)+\sum_{\gamma\in\Gsub}\phi_{-}(\gamma)\ \phi(\Gamma/\gamma)\big),\quad \text{for all }\Gamma\in\cH.
\label{eq:PhiMinus}
  \end{equation}
Let us mention here that this twisted antipode implements the BPHZ subtraction (see Section~4).
The renormalized Feynman amplitude $\phi_{+}$ of a graph $\Gamma\in\cH$ is given by
\begin{equation}
\label{rena}  
\phi_{+}(\Gamma)=\phi_{-}\ast\phi(\Gamma).
  \end{equation}
Moreover, Bogoliubov's subtraction operator is given by (see, for example, 
\cite{dirk-fields} for details)
\begin{equation}
\label{r-ck}
\bar R (\Gamma) = \phi (\Gamma) + \sum_{\gamma\in\Gsub} \phi_-(\gamma)\phi(\Gamma/\gamma).
\end{equation}
Let us now illustrate this with the example of the divergent Feynman graph of Figure~\ref{big}. In that case, the expression \eqref{r-ck} above becomes
\begin{equation}
\label{a11}
\bar R (\phi(\Gamma))=\phi (\Gamma)+\phi_-(\gamma)\phi(\Gamma/\gamma),
\end{equation}
where the only subdivergent graph $\gamma$ is given by the internal edges $1$ and $2$ (see again Figure~\ref{big}). The above definition of the 
twisted antipode leads to
\begin{equation}
\label{a21}
\phi_{-}(\gamma)=-T(\phi(\gamma))
\end{equation}
(since the subgraph $\gamma$ is primitive --- the recursion in the definition of the twisted antipode stops).
Inserting the expression \eqref{a21} for the twisted antipode in  \eqref{a11}  finally gives
\begin{equation}
\label{a33}
\bar R (\phi(\Gamma))=\phi (\Gamma)-T(\phi(\gamma))\phi(\Gamma/\gamma),
\end{equation}
which is, as announced, the same expression as \eqref{ex-qft}, the expression computed using non-algebraic formulas.  Note that the explicit form of the projection $T$  is given by the choice of the renormalization scheme. In Section~4, working in configuration space, we have used explicit Taylor expansion operators $\tau$.

\medskip

For the sake of completeness, let us recall that a combinatorial Hopf algebraic structure like the one described in this section can be also defined for more involved quantum field theories, such as gauge theories \cite{CDE3,gauge2}.

Moreover, let us also state that a slightly different combinatorial Hopf algebra --- the {\it core Hopf algebra} --- was defined in \cite{core} (and independently in \cite{thomas} for vacuum graphs).
Its coproduct does not sum over the class of {\it superficially divergent} subgraphs but over {\it all} subgraphs of the respective graph:
\begin{equation}
\Delta \Gamma = \Gamma \otimes 1 + 1 \otimes \Gamma + \sum_{\gamma\subset \Gamma} \gamma \otimes \Gamma/\gamma, \ \forall\Gamma\in\cH.
\end{equation}
The definitions of the coproduct and of the counit follow in a straightforward manner. One can then verify that this structure is also a Hopf algebra structure. The cohomology of this Hopf algebra was also investigated in detail \cite{core2A,core2B}.

This core algebra contains the renormalization Hopf algebra (presented in this section) as a quotient algebra. Moreover, it can be seen as the Connes--Kreimer Hopf algebra of a QFT model formulated in infinite dimension (see again \cite{core2A,core2B} for more details).

Let us end this section by mentioning that the study of the cohomology of the Connes--Kreimer Hopf algebras allows one to express the combinatorial Dyson--Schwinger equations as power series in some appropriate insertion operators of graphs (which, from a mathematical point of view, are Hochschild one-cocyles). The interested reader may consult the articles \cite{CDE1,CDE3,CDE2} or the thesis  \cite{teza-yeats}. For various results on these combinatorial Dyson--Schwinger equations in the Hopf algebra of decorated rooted trees, one may also consult \cite{foissy1,foissy2}.

\section{Ribbon graphs 
and scalar QFT on the non-commutative Moyal space}
\label{sec:gen}
\renewcommand{\theequation}{\thesection.\arabic{equation}}
\setcounter{equation}{0}

In this section we present the generalization of the results of the previous sections to the case of ribbon graphs (where the Tutte polynomial generalizes to the Bollob\'as--Riordan polynomial) and the case of $\Phi^4$ scalar QFT on the $4$-dimensional non-commutative Moyal space.

\subsection{Ribbon graphs --- the Bollob\'as--Riordan polynomial}

One has the following definition, generalizing in a natural manner the definition of a graph.

\begin{definition}
\label{defr}
A {\it ribbon graph} $\Gamma $ is an orientable surface with
boundary represented as the union of closed disks, also called
vertices, and  ribbons also called edges, such that:
the disks and the ribbons intersect in disjoint line
segments,
 each such line segment lies on the boundary of precisely one
disk and one ribbon and, finally, every ribbon contains two such line segments.
\end{definition}

Let us also mention that ribbon graphs can also be defined as graphs equipped with, for each vertex, a cyclic ordering of the edges incident to that vertex, or as graphs embedded in surfaces (this last definition was actually the object on which B.~Bollob\'as and O.~Riordan defined their generalization of the Tutte polynomial in \cite{br1,br2}).

The definition above can be extended such that a distinct type of edge --- the external edge --- is permitted. This type of edge contains only one line segment given by the intersection between a ribbon and a disk (see Definition~\ref{defr}).

\begin{figure}[bth]
\centerline{\includegraphics[width=2cm, angle=0]{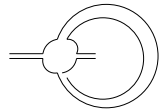}}
\caption{An example of a ribbon graph with  one vertex, one internal edge, and two external edges.
\label{ex1}}
\end{figure}

An example of such a graph is given in Figure~\ref{ex1}.

\begin{definition}
A {\it face} of a ribbon graph is a connected component of its boundary as a surface.
\end{definition}

For example, the graph of Figure~\ref{ex1} has two faces.

If we glue disks along the faces, we obtain a closed Riemann surface
whose genus is also called the genus of the graph.

\begin{definition}
The ribbon graph is called {\it planar} if this Riemann surface has genus zero. 
\end{definition}

For example, the graph of Figure~\ref{ex1} is planar. 

\begin{definition}
A planar graph is called {\it regular} if the number of faces broken by external edges is equal to $1$.
\end{definition}

The graph of Figure~\ref{ex1} is planar irregular.

\medskip

As already stated in the Introduction, ribbon graphs are also known as {\it maps}  (see for example \cite{maps} and references therein). Planar regular maps are also known as {\it outer maps}. 

\bigskip

\begin{definition}
The Bollob\'as--Riordan polynomial of a ribbon graph $G$ is defined by
\begin{equation}
R_\Gamma(x,y,z)
:=\sum_{H\subset E}(x-1)^{r(\Gamma)-r(H)}y^{n(H)}z^{k(H)-F(H)+n(H)}.
\end{equation}
\end{definition}

In the above definition, the symbol $F(H)$ stands for the number of components of the boundary of the respective subgraph $H$ (the number of faces). The supplementary variable $z$ is required to keep track of the topological information (the graph genus or the number of faces).


As the Tutte polynomial, the Bollob\'as--Riordan polynomial also obeys a de\-le\-tion/con\-trac\-tion relation. 

\begin{theorem}
\label{delcontrbollo}
Let\/ $\Gamma$ be a ribbon graph. Then
\begin{equation}
\label{cdBR}
R_\Gamma=R_{\Gamma/e}+R_{\Gamma-e}
\end{equation}
for any regular edge $e$ of $G$, 
and
\begin{equation}
R_\Gamma=x R_{\Gamma/e}
\end{equation}
for every bridge of\/ $\Gamma$.
\end{theorem}

There are further analogies to
the Tutte polynomial, in the sense that, defining the Bollob\'as--Riordan polynomial on terminal forms transforms property~\eqref{cdBR} into a definition. On these terminal forms (that is, graphs with one vertex), the polynomial is defined by
\begin{equation}
R_\Gamma (y,z):= \sum_{H \subset {\Gamma}}  y^{|E(H)|} z^{2g(H)},
\end{equation}
since, in this case, one has $k(H)-F(H)+n(H)=2g(H)$.

\medskip

A multivariate version of the Bollob\'as--Riordan polynomial exists as well in the literature, namely
\begin{equation}
Z_\Gamma(x,\beta,z)
=\sum_{H\subset E}x^{k(H)}\left(\prod_{e\in H}\beta_e\right)z^{F(H)}.
\end{equation}

This version also satisfies a deletion/contraction relation.

\medskip

Finally, let us mention that a signed version of the Bollob\'as--Riordan 
was defined in
 \cite{signedBR}; this is a three-variable polynomial defined on {\it signed ribbon graphs}, that is, ribbon graphs for which an element of the set $\{+,-\}$ is assigned to each edge.
A partial duality with respect to a spanning subgraph was also defined \cite{pd}; this allows one to prove that the Kauffman bracket of a virtual link diagram is equal to the signed Bollob\'as--Riordan polynomial of some ribbon graph constructed from a state of the respective virtual link diagram (the interested reader is referred to  \cite{pd} for details on this topic).
Moreover, the properties of the multivariate version of this signed  Bollob\'as--Riordan polynomial were analyzed in \cite{fab-br} (namely, there its invariance under the partial duality of \cite{pd} was proved). 
Finally, $4$-variable generalizations of the Bollob\'as--Riordan polynomial for ribbon graphs were defined in \cite{sur} and \cite{pol2}.

\subsection{Scalar $\Phi^{\star\, 4}$ QFT on the non-commutative Moyal space;  parametric representation}

Let us start this subsection with some explanations on the interest of mathematical physicists in QFT on non-commutative spaces.

Non-commutative QFT can be seen as a potential candidate for new physics beyond the celebrated Standard Model of elementary particle physics. Moreover, non-commutativity of space-time is believed to be a possible framework for the quantification of gravity. Thus, it was argued in \cite{dfr} that the space-time becomes non-commutative at the Planck scale (the physical energy scale where gravitational effects have the same order of magnitude as quantum effects; thus, the energy scale where gravity and quantum physics cannot ``ignore themselves" anymore).
Moreover, non-commutative QFT has been shown to be related to loop quantum gravity and to group field theory \cite{fln3,fln1,fln2}, some of the major candidates for a quantified theory of gravity (more details on group field theory are also given in the last section of this survey). 
Non-commutative models have also been shown to be effective theories of certain models  of string theories and of matrix theories \cite{csw1,csw2}. Finally, changing point of view, it was argued that non-commutative models can be better adapted for describing the physics of effective non-local interactions, as it is the case in the fractional quantum Hall effect \cite{poly1,poly2}. For more details on non-commutative QFT, the interested reader may consult the review paper \cite{rev-riv}.

\medskip

From a mathematical point of view, the Moyal space, one of the simplest cases of such a non-commutative setting, can be defined in the following way.

\begin{definition}
The {\it Moyal algebra} is the linear space of smooth and rapidly decreasing functions ${\mathcal S} ({\mathbb R}^D)$ equipped with the 
 {\it Moyal product}
\begin{equation}
 (f\star g)(x)=\int \frac{d^D k}{(2\pi)^D} d^Dy \, f(x+\frac 12 \Theta \cdot k) g(x+y) e^{ik\cdot y}, \forall x\in \RR^D,
\end{equation}
where $\Theta$ is a $D$-dimensional skew-symmetric matrix.
\end{definition}

\begin{remark}
\noindent
\begin{enumerate}
\item The Moyal product is non-local, non-commutative, but it is still associative.
\item Taking the matrix $\Theta$ to be vanishing leads to the usual commutative multiplication of functions. 
\end{enumerate}
\end{remark}

The Moyal algebra defined above can be extended by duality, considering the product of a tempered distribution with a Schwartz-class function. The identity, the polynomials, the $\delta$ function, and its derivatives can then be shown to belong to this extended Moyal algebra (see \cite{gbv} for details).
One can then prove that
\begin{equation}
[x_\mu, x_\nu]_\star:=x_\mu\star x_\nu-x_\nu\star x_\mu=i \Theta^{\mu \nu},\quad \text{for all }\mu,\nu =1,2,\ldots,D.
\end{equation}

Let us note that this relation can also be used as a definition of the Moyal space.

In this survey we take $D=4$ (the usual $4$-dimensional space-time). The {\it non-commutativity matrix} $\Theta$ is taken to be
$$
\Theta = \begin{pmatrix} \Theta_2 & 0 \\ 0 & \Theta_2 \end{pmatrix},\ \
\Theta_2 = \begin{pmatrix} 0 & -\theta \\ \theta & 0 \end{pmatrix}. 
$$

Non-commutative quantum field theories are then implemented by replacing the usual commutative multiplication of fields by the non-commutative $\star$-product. Thus, the action of the non-commutative $\Phi^{\star\, 4}$ model reads
\begin{equation}
\label{act-nc}
S[\Phi]=\int_{\RR^4} d^4 x \left[ 
\frac{1}{2} \sum_{\mu=1}^4\left(\partial_\mu  \Phi\star \partial_\mu\Phi\right) (x)
 +\frac{1}{2}m^2{\Phi}^{\star\, 2} (x)+ V_{\text{int}}^\star[\Phi]
\right],
\end{equation} 
where the non-commutative interaction potential is given by
\begin{equation}
V_{\text{int}}^\star[\Phi]=\frac{\lambda}{4} \Phi^{\star\, 4} (x).
\end{equation}

As in the commutative case, one can apply the Fourier transform in order to
transfer this action to one on momentum space.
A first consequence of the use of the Moyal product is that its use does not change the propagation part of the action. This is a direct consequence of the following mathematical property:
\begin{equation}
\int d^4 x \, (\Phi  \star \Psi) (x) = \int d^4 x \, \Phi (x) \, \Psi (x), 
\quad \text{for all } \Phi\text{ and } \Psi.
\end{equation}
Thus, the propagator of the theory has the same form as in the commutative case; this is particularly useful for implementing the parametric representation in a similar way.

Nevertheless, the use of the Moyal product at the level of the interaction terms changes completely the situation: the non-commutative interaction is not equivalent to the commutative one (as this is also the case for the propagation).
The interaction part
no longer preserves the invariance under permutation of external fields. This invariance is 
restricted to {\it cyclic permutations} only. Furthermore, there exists a basis --- the matrix base --- of the Moyal algebra where the Moyal product takes the form of an ordinary  matrix product. For these reasons
the associated Feynman graphs are {\it ribbon graphs}, just like the ones described in the previous subsection. 

\medskip

For more details on the implementation of QFTs on the non-commutative Moyal space, the interested reader is referred to \cite{rev-riv} and references therein.

\bigskip

Let us now move on to the  implementation of the {\it non-commutative parametric representation}. This can be achieved
using the same trick \eqref{trick} as in the commutative case. 
The difference comes however from the fact that the interaction part is non-trivially dependent on the non-commutativity parameter $\theta$, and this dependence will consequently appear in the form of the two polynomials which give the parametric representation of the respective Feynman integral (we denote these polynomials by $U^\star (\alpha, \theta)$ and $V^\star(p_{{\rm ext}},\alpha, \theta)$).

Using a Grassmannian development of the Pfaffians obtained through these Gaussian integrations, one can prove (see again \cite{io-BR}) that the Feynman integral of a generic graph $\Gamma$ reads
\begin{equation} 
\phi(\Gamma)=
\int_0^{\infty} 
\frac{e^{- V^\star(p_{{\rm ext}},\alpha,\theta)/U^\star (\alpha) }}{U^\star (\alpha,\theta)^{2}} 
\prod_{\ell=1}^{|E|}  
( e^{-m^2 \alpha_\ell} 
d\alpha_\ell 
).
\end{equation}

In order to present the exact form of the polynomials  $U^\star_ \Gamma(\alpha, \theta)$ and $V^\star_\Gamma(p_{{\rm ext}},\alpha, \theta)$,we need the following definition (see again \cite{io-BR} for details).

\begin{definition}
\noindent
\begin{enumerate}
\item A $\star$-tree of a connected ribbon graph  is a sub-ribbon graph with  one boundary component.
\item A two $\star$-tree is a sub-ribbon graph with two boundary components.
\end{enumerate}
\end{definition}

Note that these definitions are non-trivial generalizations of the notions of trees and two-trees, respectively (see Definition~\ref{toate}(8) and (9)). Let us also mention that $\star$-trees are also known as {\it quasi-trees} in the knot theory literature and as unicellular maps in the combinatorics literature (see for example \cite{fab-br} and \cite{fab-br2} and references therein).

We have
\begin{equation}
\label{min1}
U^\star_\Gamma(\alpha, \theta)=\left(\frac{\theta}{2}\right)^{b} \sum_{{\mathcal T}^\star\; \star\text{-}{\mathrm tree}} \prod_{e\notin {\mathcal T}^\star} 2 \frac{\alpha_e}{\theta},
\end{equation}
where 
$b=F-1+2g$. 

The case of the polynomial  $V^\star(p_{{\rm ext}},\alpha,\theta)$ is more involved, because the non-com\-mu\-ta\-tive setting presented here allows for it to be complex. We have the following theorems.

\begin{theorem}[{\sc\cite[Theorem 5.2]{io-BR}}]
\label{thr}
The real part of the  polynomial  $V^\star_\Gamma(p_{{\rm ext}},\alpha,\theta)$ reads
\begin{equation}
\label{vr}
{\mathcal X}^{\star}_\Gamma(p_{{\rm ext}},\alpha,\theta)
=\left(\frac{\theta}{2}\right)^{b+1}\sum_{{\mathcal T}^\star_2\; \text{\em two-}\star{\text{\em-tree}}}\prod_{e\notin {\mathcal T}_2^\star}2\frac{\alpha_e}{\theta}(p_{{\mathcal T}^\star_2})^2,
\end{equation}
where $p_{{\mathcal T}^\star_2}$ is the sum of 
the external momenta entering one of the two faces of the two-$\star$-tree ${\mathcal T}_2^\star$.
\end{theorem}

Note that we say that an external momentum enters a face
 if the respective external edge is hooked to one of the vertices of the face (we also consider
in the sum above, for each such external momentum, 
 the sign $+$ if the momentum is ingoing, respectively $-$ 
if it is outgoing).  
As in the commutative case,
the choice of the face in the above theorem is irrelevant (by momentum conservation).

\begin{theorem}[{\sc\cite[Theorem 5.3]{io-BR}}]
\label{thi}
The imaginary part of the polynomial $V^\star_\Gamma(p_{{\rm ext}},\alpha,\theta)$ reads
\begin{equation}
\label{vim}
{\mathcal Y}^{\star}_\Gamma(p_{{\rm ext}},\alpha,\theta)=\left(\frac{\theta}{2}\right)^{b}\sum_{{\mathcal T}^\star\; \star{\text{\em-tree}}}\prod_{e\notin {\mathcal T}^\star}2\frac{\alpha_e}{\theta} \psi(p),
\end{equation}
where $\psi(p)$ is the phase obtained by following 
the momenta entering the face of the $\star$-tree ${\mathcal T}^\star$.
\end{theorem}

Let us give some explanations on the definition of  $\psi(p)$. 
We denote the (unique) face of the $\star$-tree by $f$. 
One can then define an incidence matrix $(\tilde\epsilon_{f,i})$ between this  face and the external edges of the graph (in the same way one defines the usual incidence matrix between vertices and edges).  Following 
the external momenta entering the face, one can define an order relation on these momenta.
One then has $\psi(p)=\sum_{i<j} (\tilde \epsilon_{f,i}p_{i})\wedge (\tilde \epsilon_{f,i}p_{j})$, 
where  $p_i\wedge p_j=\frac{1}{\theta}p_i^{\mu}\Theta_{\mu\nu}p_j^{\nu}$. 
This expression generalizes in a straightforward manner 
the expression which gives, as a function of the incoming/outgoing momenta, 
the phase corresponding to a Moyal vertex.

\subsection{Relation between the multivariate Bollob\'as--Riordan
polynomial and the polynomials of the parametric representation}

It can be shown that also
the polynomials $U^\star_\Gamma(\alpha,\theta)_\Gamma$ and  $V^\star_\Gamma(p_{{\rm ext}},\alpha,\theta)$ satisfy a deletion/contraction relation. 
More precisely, we have the following result.

\begin{theorem}[{\sc\cite[Theorem 5.1]{io-BR}}]
\label{grasstheo2}
For any semi-regular edge $e$, we have
\begin{equation} 
\label{delcontr3}
U^\star_\Gamma(\alpha,\theta) = \alpha_e U^\star_{\Gamma- e}(\alpha,\theta) + U^\star_{\Gamma/ e}(\alpha,\theta).
\end{equation}
\end{theorem}

One can moreover prove the following relation between the multivariate Bollob\'as--Riordan polynomial and the polynomial $U^\star_\Gamma(\alpha,\theta)$ of the parametric representation:
\begin{equation}
\label{limita-nc}
U^{\star}_{\Gamma}(\alpha,\theta)=(\theta/2)^{\frac 12 (|V|-|E|-1)}\Big(\prod_{e\in E}\alpha_{e}\Big)\times
\lim_{w\rightarrow0}w^{-1}Z_{\Gamma}\Big(1, {\textstyle\frac{\theta}{2\alpha_{e}}},w\Big).
\end{equation}
Note that the commutative limit ($\theta\to 0$) gives the usual $\Phi^4$ model on the commutative $\RR^4$. This can be seen as well in the parametric representation, where $U^{\star}_{\Gamma}(\alpha,\theta)$ is given by  $U^{}_{\Gamma}(\alpha)$ plus some noncommutative corrections which factor $\theta^2$ (corrections thus vanishing in the commutative limit). This directly applies to equation~\eqref{limita-nc} (where one obviously has to take the limits in the appropriate order on the right-hand side), and it leads to
\begin{equation}
U^{}_{\Gamma}(\alpha)=\lim_{\theta\rightarrow 0}U^{\star}_{\Gamma}(\alpha,\theta)=\lim_{\theta\rightarrow 0}\left((\theta/2)^{\frac 12 (|V|-|E|-1)}\Big(\prod_{e\in E}\alpha_{e}\Big)\times
\lim_{w\rightarrow0}w^{-1}Z_{\Gamma}\Big(1, {\textstyle\frac{\theta}{2\alpha_{e}}},w\Big)\right).
\end{equation}

 The second polynomial, $V^\star_\Gamma(\alpha,p_{{\rm ext}},\theta)$,  can also be related to the multivariate Bollob\'as--Riordan polynomial in a similar way (see again \cite{io-BR} for details).

\medskip

For the sake of completeness, let us also mention that the parametric representation of a 
non-commutative renormalizable scalar model, 
different from \eqref{act-nc} (the latter being actually non-renormalizable)
 --- the Grosse--Wulkenhaar model \cite{GW} --- was fully implemented. 
Its polynomials do not obey a deletion/contraction relation, but they obey some other types of recursive relations  \cite{pol2}. 
Algebraic topological properties of the parametric representation of the Grosse--Wulkenhaar model (as well as of the parametric representation of the commutative model \eqref{act}) were studied in \cite{mar}.

\section{Combinatorial Connes--Kreimer Hopf algebras for renormalizable 
$\Phi^{\star\, 4}$ QFT on the non-commutative Moyal space}
\label{sec:ck-fab}
\renewcommand{\theequation}{\thesection.\arabic{equation}}
\setcounter{equation}{0}




The definition of the Hopf algebra of non-commutative Feynman graphs which 
drives the combinatorics of  renormalization is formally the same as in the commutative case.  
Nevertheless, in the definition of the coproduct one has to sum
over
the superficially divergent subgraphs of the respective graph, and
this class of superficially divergent graphs takes now explicitly
into account the topology of the graph.
Thus, using various QFT techniques, it was proved in \cite{GW} 
that the superficially divergent graphs 
of the renormalizable Grosse--Wulkenhaar 
non-commutative model are the {\it planar regular} ribbon graphs with two and four external edges.

\medskip

We denote by $\cH^\star$ the unital, associative algebra freely generated by $1$PI ribbon 
Feynman graphs of the 
renormalizable Grosse--Wulkenhaar 
$\Phi^{\star\, 4}$ model  (including the empty set, which we denote by $1_{\cH^\star}$). The product $m$ is bilinear, commutative and given by the operation of disjoint union. Let the coproduct $\Delta:\cH^\star\to\cH^\star\otimes\cH^\star$ be defined by
\begin{equation}
\Delta \Gamma = \Gamma \otimes 1_{\cH^\star} + 1_{\cH^\star}\otimes \Gamma + \sum_{\gamma\in \Gsub} \gamma \otimes \Gamma/\gamma, \quad \text{for all }\Gamma\in\cH^\star,
\label{eq:NCcoproduct2}
\end{equation}
where, as earlier, we write $\Gsub$ for the set of two and four external edges planar regular 
sub-ribbon graphs of $\Gamma$.
Furthermore, we define the counit $\veps:\cH^\star\to\KK$ by
\begin{equation}
\veps (1_{\cH^\star}) =1,\ \veps (\Gamma)=0,\ \forall \Gamma\ne1_{\cH^\star}.
\end{equation}
Finally, the antipode is given recursively by
\begin{align}
  S:\cH^\star\to&\cH^\star\label{eq:Antipode2}\\
  \Gamma\mapsto&-\Gamma-\sum_{\gamma\in\Gsub}S(\gamma)\Gamma/\gamma.\nonumber
\end{align}

Then the following result holds true.

\begin{theorem}[{\sc\cite[Theorem 4.1]{fab}}]
\label{thm:hopf-algebra-nonComm}
The quadruple $(\cH^\star,\Delta,\veps,S)$ is a Hopf algebra.
\end{theorem}

\begin{figure}[bth]
\centerline{\includegraphics[width=5cm]{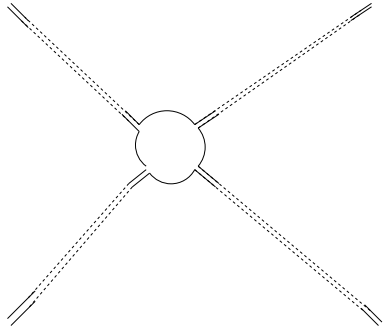}}
\caption{Insertion of a regular ribbon graph with four external edges. The cyclic ordering is respected by the chosen gluing data.
\label{fig:insertion}}
\end{figure}

The main difficulty in the proof of this theorem 
(in comparison with the original Connes--Kreimer theorem) 
comes from the necessity of defining gluing data 
which respect to the ribbon graph cyclic ordering. 
Figure~\ref{fig:insertion} shows an example of an insertion 
of a graph with four external edges with gluing data 
which respects the cyclic ordering 
(the dotted lines showing how the four external legs 
of the subgraph to be inserted and the four edges of the chosen vertex 
are put together --- the insertion place).

\bigskip

The renormalized Feynman amplitude is defined in analogy to the commutative case (see equation~\eqref{rena}).

\medskip

Let us end this section by recalling that in \cite{io-kreimer} 
the role played by Hochschild cohomology of 
these types of Hopf algebras in non-commutative QFT 
was studied in detail; non-trivial 
examples of graphs with one or two independent cycles have been explicitly worked out.

\section{Perspectives --- combinatorics of quantum gravity tensor models}

One of the main perspectives of the combinatorial approaches presented here is their extension for the study of the  properties of quantum gravity models. 
The group field theory formalism of quantum gravity (for general reviews,
see \cite{gft1,gft2,gft3,gft4}) is the most adapted for such a study, since it is formulated as a QFT. These models were developed as a generalization of
two-dimensional matrix models (which naturally make use of ribbon graphs just as non-commutative quantum field theories do) to the three- and four-dimensional cases.

\begin{figure}
\begin{center}
\includegraphics[scale=0.3,angle=0]{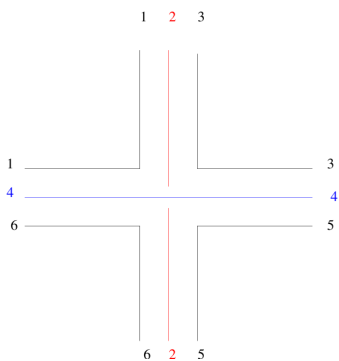}
\caption{A possible vertex of a three-dimensional quantum gravity tensor model.}
\label{fig:3d}
\end{center}
\end{figure}

\begin{figure}
\begin{center}
\includegraphics[scale=0.3,angle=0]{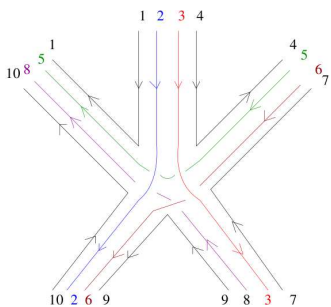}
\caption{A possible vertex of a four-dimensional quantum gravity tensor model.}
\label{fig:4d}
\end{center}
\end{figure}

The natural candidates for generalizations of matrix models in higher dimensions ($>2$)
are {\it tensor models}. In the combinatorially simplest case, the elementary cells that, by gluing
together form the space itself, are the $D$-simplices ($D$ being the dimension of space). Since a
$D$-simplex has $D+1$ facets on its boundary, the backbone of group field theoretical models in
dimension $D$ should be some abstract $\Phi^{D+1}$ interaction on rank $D$ tensor fields $\phi$ (see Figures~\ref{fig:3d} and \ref{fig:4d}, which represent, for $D=3$ and $D=4$, the vertex of these Feynman quantum gravity tensor graphs).

In \cite{gp1,gp2}, some propositions for generalizations of the Bollob\'as--Riordan polynomial to the level of these rank three tensor models have been established, and it was proven that they respect a  deletion/contraction property. 
Due to the increased topological complexity of the graphs, a different way of contracting edges had to be used.

Furthermore, a polynomial for triangulations (no graphs being taken into considerations) has also been recently proposed in \cite{surpriza}. It is thus interesting to understand if some connection between these propositions exists and if any of them could be related to some parametric representation of these quantum gravity models, as shown here to be the case for both commutative and non-commutative QFTs.

Let us also stress the fact that different topological and analytic
insights in this type of formulation of quantum gravity models have
been investigated in the recent literature (see
\cite{ultimele7,fln3,ultimele3,ultimele5,ultimele1,ultimele6,ultimele2,ultimele8} and references therein).

\medskip

It would be interesting to investigate whether or not these models can be related to the operation of cabling in knot theory and thereby to the colored Jones polynomials. Within the quantum gravity framework described in this section, (edge-)colored tensor models have been recently proposed and investigated (see \cite{ceaslov} and references therein). On the other hand, the relation of parallel cabling to combinatorics in general and to the Bollob\'as--Riordan polynomial of a ribbon graph in particular was analyzed in \cite{ultim}.

\bigskip

For the sake of completeness, let us end this survey by mentioning
other types of mathematical perspectives, such as number theoretical
conjectures related to Feynman integrals, definition of knot
invariants, etc. The interested reader may for example consult
\cite{altii2} and \cite{b-final3,b-final4,b-final2,b-final1} (and
references therein) for various other developments related to the
combinatorics of graph polynomials and QFT. 

\bigskip

\section*{Acknowledgments}
 The author acknowledges the grants CNRS PEPS ``CombGraph'',  PN 09 37 01 02 and  CNCSIS ``Tinere echipe'' 77/04.08.2010.

\appendix

\section{Grassmann variables, determinants and Pfaffians}
\label{grass}

Grassmann variables $\chi_1, \chi_2,\dots, \chi_n$
are defined through their anticommutation relations
\begin{equation}
\label{defg}
 \chi_i \chi_j = - \chi_j \chi_i, \quad \text{for all }  i\text{ and }  j.
\end{equation}
A direct consequence of this relation is that any function of these variables is a polynomial
with highest degree one in each variable. 
The rules of Grassmann integrations are then
\begin{equation}
\label{intg}
  \int d\chi = 0\quad \text{and}  \quad \int \chi\, d\chi = 1.
\end{equation}

The determinant of any $n$-dimensional square matrix can be expressed as
a Grassmann Gaussian integral over $2n$ independent
Grassmann variables which one can denote by 
$\bar \psi_1, \bar\psi_2,\ldots , \bar \psi_n$, $\psi_1,\psi_2, \ldots ,  \psi_n$.
 (The reader should be warned that the bars have nothing to do with complex conjugation.) More precisely, we have
\begin{equation} 
\det M = \int \prod d\bar \psi_i d\psi_i   e^{-\sum_{ij} \bar\psi_i M_{ij} \psi_j   } .
\end{equation}

The Pfaffian $\mathrm{Pf} (A)$ of a \emph{skew-symmetric} 
matrix $A$ is defined by
\begin{equation}
\det A = [\mathrm{Pf} (A)]^2 ,
\end{equation}
where the sign of the term $A_{1,2}A_{3,4}\cdots$ is fixed to be $+1$.

\begin{proposition}[{\sc\cite[Proposition 3]{malek}}]
The Pfaffian of an $n$-dimensional skew-sym\-metric matrix $A$ is given by
\begin{equation}
\mathrm{Pf} (A) =\int \prod d\chi_i
e^{-\sum_{i<j}\chi_i A _{ij}\chi_j}
= \int d\chi_1\,d\chi_2\cdots d\chi_n e^{-\frac{1}{2}\sum_{i,j}\chi_i A _{ij}\chi_j} .
\label{pfaff}
\end{equation}
\end{proposition}

\begin{lemma}[{\sc\cite[Lemma 2.1]{io-BR}}]
\label{lemma2}
The determinant of a matrix $D+A$, where $D$ is
diagonal and $A$ skew-symmetric, can be written in the form
\begin{equation}  
\det (D+A) = \int \prod d\chi_i d \omega_i e^{-\sum_i  \chi_i D_{ii} \omega_i - \sum_{i <j}
\chi_i A_{ij} \chi_j + \sum_{i <j}  \omega_i A_{ij} \omega_j } .
\end{equation}
\end{lemma}


\begin{thebibliography}{99}

\bibitem {malek}
A. Abdesselam,
{\em The {G}rassmann--{B}erezin calculus and theorems of the matrix-tree type},
Adv. Appl. Math. {\bf 33} (2004), 51--70.

\bibitem{mar}
  P.~Aluffi and M.~Marcolli,
{\em Feynman motives of banana graphs}, 
 Commun. Number Theory Phys. {\bf 3} (2009), 1--57;
  {\tt ar$\chi$iv:0807.1690}.

\bibitem{fln3}
  A.~Baratin and D.~Oriti,
{\em Group field theory with non-commutative metric variables},
  Phys.\ Rev.\ Lett.\  {\bf 105} (2010), 221302, 4~pp.;
  {\tt ar$\chi$iv:1002.4723}.

\bibitem{ultimele7}
J.~Ben Geloun, R.~Gurau, V.~Rivasseau,
{\em EPRL/FK group field theory},
Europhys. Lett. {\bf 92} (2010), 60008, 6~pp.;
  {\tt ar$\chi$iv:1008.0354}.

\bibitem{ultimele3}
 J.~Ben~Geloun, T.~Krajewski, J.~Magnen and V.~Rivasseau,
  {\em Linearized Group Field Theory and Power Counting Theorems},
  Class.\ Quant.\ Grav.\  {\bf 27}  (2010), 155012, 14~pp.;
  {\tt ar$\chi$iv:1002.3592}.


\bibitem{CDE1}
  C.~Bergbauer and D.~Kreimer,
{\em Hopf algebras in renormalization theory: locality and Dyson--Schwinger
equations from Hochschild cohomology},
  IRMA Lect.\ Math.\ Theor.\ Phys.\  {\bf 10} (2006), 133--164;
  {\tt ar$\chi$iv:hep-th/0506190}.

\bibitem{altii1}
  S.~Bloch, H.~Esnault and D.~Kreimer,
  {\em On motives associated to graph polynomials},
  Commun.\ Math.\ Phys.\  {\bf 267} (2006), 181--225;
  {\tt ar$\chi$iv:math/0510011}.

\bibitem{core}
S. Bloch and D. Kreimer, {\em Mixed Hodge structures and renormalization in physics}, 
Commun. Number Theory Phys. {\bf 2} (2008), 
637--718; 
{\tt ar$\chi$iv:0804.4399}.

\bibitem{br1}
B. Bollob\'as and O. Riordan, {\em A polynomial invariant of graphs on
  orientable surfaces}, Proc. London Math. Soc. {\bf 83} (2001), 513--531.

\bibitem{br2}
B. Bollob\'as and O. Riordan, 
{\it A polynomial of graphs on surfaces},
Math. Ann. {\bf 323} (2002), 81--96.

\bibitem{ultimele5}
V.~Bonzom, M.~Smerlak,
{\em Bubble divergences from cellular cohomology},
  Lett.\ Math.\ Phys.\  {\bf 93} (2010), 295--305;
  {\tt ar$\chi$iv:1004.5196}.


\bibitem{b-final3}
F. Brown, {\em On the periods of some Feynman integrals}, 
preprint;
{\tt ar$\chi$iv:0910.0114}.


\bibitem{altii2}
  F.~Brown,
{\em The massless higher-loop two-point function},
  Commun.\ Math.\ Phys.\  {\bf 287} (2009), 925--958;
  {\tt ar$\chi$iv:0804.1660}.

\bibitem{b-final4}
F.~Brown and D.~Kreimer,
{\em Angles, scales and parametric renormalization}, preprint;
  {\tt ar$\chi$iv:1112.1180}.

\bibitem{b-final2}
F. Brown and O. Schnetz, {\em A K3 in $\phi^4$}, 
Duke Math. J. (to appear);
 {\tt ar$\chi$iv:1006.4064}.

\bibitem{b-final1}
 F.~Brown and K.~A. Yeats,
{\em Spanning forest polynomials and the transcendental weight of Feynman graphs},
  Commun.\ Math.\ Phys.\  {\bf 301} (2011), 357--382;
  {\tt ar$\chi$iv:0910.5429}.



\bibitem{fab-br2}
A. Champanerkar, I. Kofman and N. Stoltzfus, {\em Quasi-tree expansion
  for the Bollob\'as--Riordan-Tutte polynomial}, Bull. London
Math. Soc. {\bf 43} 2011), 972--984.


\bibitem{maps}
G. Chapuy, {\em Combinatoire bijective des cartes de genre sup\'erieur}, Ph.D. thesis (in French), \'Ecole Polytechnique, Paris, 2009;
available at\newline
{\tt http://www.liafa.univ-paris-diderot.fr//\lower0.5ex\hbox{\~{}}chapuy/articles/these-chapuy.pdf}.


\bibitem{pd}
S. Chmutov, {\em Generalized duality for graphs on surfaces and the
  signed Bollob\'as--Riordan polynomial}, J. Combin.
  Theory Ser.~B, {\bf 99} (2009), 617--638, 
{\tt ar$\chi$iv:0711.3490}.

\bibitem{signedBR}
S. Chmutov and I. Pak, {\em The Kauffman bracket of virtual links and
  the Bollob\'as--Riordan polynomial}, Moscow Math. J. {\bf 7}
(2007), 409--418.


\bibitem{ck}
  A.~Connes and D.~Kreimer,
 {\em Renormalization in quantum field theory and the Riemann--Hilbert  problem,
 I: the Hopf algebra structure of graphs and the main theorem},
  Commun.\ Math.\ Phys.\  {\bf 210} (2000), 249--273;
  {\tt ar$\chi$iv:hep-th/9912092}.


\bibitem{csw1}
  A.~Connes, M.~R.~Douglas and A.~S.~Schwarz,
  {\em Noncommutative geometry and matrix theory: compactification on tori},
  J. High Energy Phys. {\bf 9802} (1998), Paper~3, 35~pp.;
  {\tt ar$\chi$iv:hep-th/9711162}.

\bibitem{tutte2}
H. Crapo, {\em The Tutte polynomial}, aequationes math. {\bf 3}
(1969), 211--229.

\bibitem{Dascalescu}
S. D\u asc\u alescu, C. N\u ast\u asescu and S. R\u aianu, {\em Hopf Algebras.
An Introduction}, Monographs and Textbooks in
Pure and Applied Mathematics, vol.~235, Marcel Dekker, Inc., 
New York, 2001.

\bibitem{dfr}
  S.~Doplicher, K.~Fredenhagen and J.~E.~Roberts,
{\em Space-time quantization induced by classical gravity},
  Phys.\ Lett.~B {\bf 331} (1994), 39--44.

\bibitem{ger}
 G.~H.~E.~Duchamp, P.~Blasiak, A.~Horzela,  K. A. Penson 
and A. I. Solomon,
 {\em Hopf algebras in general and in combinatorial physics: 
a practical introduction}, unpublished manuscript;
    {\tt ar$\chi$iv:0802.0249}.

\bibitem{patras}
K.~Ebrahimi-Fard and F.~Patras,
  {\em Exponential renormalization}, Ann. Henri Poincar\'e
{\bf 11} (2010), 943--971;
 {\tt ar$\chi$iv:1003.1679}.

\bibitem{j0A}
J. Ellis-Monaghan and C. Merino, {\em Graph polynomial and their
  applications. I. The Tutte polynomial}, 
in: Structural Analysis of Complex Networks, Matthias
  Dehmer (ed.), Birkh\"auser/Springer, New York, 2010, pp.~219--255;
{\tt ar$\chi$iv:0803.3079}. 

\bibitem{j0B}
J. Ellis-Monaghan and C. Merino,
{\em Graph polynomial and their applications. II. Interrelations and
  interpretations}, in:  Structural
Analysis of Complex Networks, Matthias Dehmer (ed.), Birkh\"auser/Springer, New York, 2010, pp.~257--292;
{\tt ar$\chi$iv:0806.4699}. 


\bibitem{foissy1}
L. Foissy, {\em Classification of systems of Dyson--Schwinger
  equations in the Hopf algebra of decorated rooted trees}, 
  Adv. Math. {\bf 224} (2010), 2094--2150. 

\bibitem{foissy2}
L. Foissy, {\em Lie algebras associated to systems of Dyson--Schwinger
  equations}, Adv. Math. {\bf 226} (2011), 4702--4730.

\bibitem{gft1}
 L.~Freidel,
{\em Group field theory: an overview},
  Int.\ J.\ Theor.\ Phys.\  {\bf 44} (2005), 1769--1783;
  {\tt ar$\chi$iv:hep-th/0505016}.

\bibitem{ultimele1}
L. Freidel, R. Gurau and D. Oriti, 
{\em Group field theory renormalization --- the 3d case: power counting of divergences}, 
Phys. Rev.~D {\bf 80} (2009), 044007, 25~pp.

\bibitem{fln1}
 L.~Freidel and E.~R.~Livine,
  {\em Effective 3-D quantum gravity and non-commutative quantum field theory},
  Phys.\ Rev.\ Lett.\  {\bf 96} (2006), 221301, 4~pp.;
  {\tt ar$\chi$iv:hep-th/0512113}.

\bibitem{gbv} 
  J.~M.~Gracia-Bond\'\i a and J.~C.~Varilly,
{\em Algebras of distributions suitable for phase space quantum
  mechanics, I},
  J.\ Math.\ Phys.\  {\bf 29} (1988), 869--879.

\bibitem{GW}
H. Grosse and R. Wulkenhaar, 
{\em Renormalization of $\phi^4$-theory on noncommutative ${{\mathbb {R}}}^4$ in the matrix base}, 
Commun. Math. Phys. {\bf 256} (2005), 305--374.



\bibitem{gp1}
R.~Gurau,
 {\em Topological Graph Polynomials in Colored Group Field Theory},
  Ann. Henri Poincar\'e {\bf 11} (2010), 565--584;
  {\tt ar$\chi$iv:0911.1945}.

\bibitem{ceaslov}
R.~Gurau and J.~P.~Ryan,
{\em Colored Tensor Models --- a review}, 
 SIGMA {\bf 8} (2012), Paper~020, 78~pp.;
  {\tt ar$\chi$iv:1109.4812}.

\bibitem{ultim}
S. Huggett, I. Moffatt and N. Virdee, {\em On the Seifert graphs of a
  link diagram and its parallels}, 
Math. Proc. Cambridge Phil. Soc. {\bf 153} (2012), 123--145;
{\tt ar$\chi$iv:1106.4197}. 

\bibitem{iz}
C. Itzykson and J.-B. Zuber, 
{\em Quantum Field Theory}, Dover Publications Inc., 2006.

\bibitem{Kassel}
C. Kassel, {\em Quantum Groups}, 
Graduate Texts in Mathematics, vol.~155,
Springer--Verlag, New York, 1995.

\bibitem{ultimele6}
T.~Krajewski, J.~Magnen, V.~Rivasseau, A. Tanasa and 
P. Vitale,
{\em Quantum corrections in the group field theory formulation of the 
Engle--Pereira--Rovelli--Livine and Freidel--Krasnov models}, 
Phys. Rev.~D {\bf 82} (2010), 124069, 20~pp.;
  {\tt ar$\chi$iv:1007.3150}.

\bibitem{thomas}
T. Krajewski and P. Martinetti, {\em Wilsonian renormalization, 
differential equations and Hopf algebras}, 
Contemp. Math. {\bf 539} (2011), 187--236; 
{\tt ar$\chi$iv:0806.4309}.

\bibitem{io-BR}
  T.~Krajewski, V.~Rivasseau, A.~Tanasa and Z.~Wang,
  {\em Topological graph polynomials and quantum field theory,
Part~I: heat kernel
  theories}, J. Noncomm. Geom. {\bf 4} (2010), 29--82;
  {\tt ar$\chi$iv:0811.0186}.

\bibitem{pol2}
T.~Krajewski, V.~Rivasseau and F.~Vignes-Tourneret,
  {\em Topological graph polynomials and quantum field theory, Part II: Mehler kernel theories},   Ann. Henri Poincar\'e 
{\bf 12} (2011), 483--545;
{\tt ar$\chi$iv:0912.5438}.


\bibitem{CDE3}
D.~Kreimer,
  {\em Anatomy of a gauge theory},
  Annals Phys.\  {\bf 321} (2006), 2757--2781;
  {\tt ar$\chi$iv:hep-th/0509135}.

\bibitem{dirk-fields}
  D.~Kreimer,
  {\em Dyson--Schwinger equations: from Hopf algebras to number theory},
  Fields Inst.\ Commun.\  {\bf 50} (2007), 225--248;
  {\tt ar$\chi$iv:hep-th/0609004}.


\bibitem{core2A}
D. Kreimer, {\em The core Hopf algebra}, in: Quanta of maths,
\'E.~Blanchard et al.\ (eds.), 
Clay Math. Proc., vol.~11, Amer. Math. Soc., Providence, 
R.I., 2010, pp.~313--321; 
{\tt ar$\chi$iv:0902.1223}. 


\bibitem{core2B}
D.~Kreimer and W.~D.~van Suijlekom,
  {\em Recursive relations in the core Hopf algebra},
Nucl. Phys.~B {\bf 820} (2009), 682--693;
 {\tt ar$\chi$iv:0903.2849}.

\bibitem{CDE2}
D.~Kreimer and K.~A. Yeats,
{\em Recursion and growth estimates in renormalizable quantum field theory},
  Commun.\ Math.\ Phys.\  {\bf 279} (2008), 401--427;
  {\tt ar$\chi$iv:hep-th/0612179}.



\bibitem{sur}
V. Krushkal, {\em Graphs, links  and duality on surfaces}, 
Combin. Probab. Comput. {\bf 20} (2011), 267--287.

\bibitem{surpriza}
V. Krushkal and D. Renardy, {\em A polynomial invariant and duality
  for triangulations}, preprint;
{\tt ar$\chi$iv:1012.1310}.

\bibitem{ultimele2}
J.~Magnen, K.~Noui, V.~Rivasseau and M.~Smerlak,
 {\em Scaling behaviour of three-dimensional group field theory},
  Class.\ Quant.\ Grav.\  {\bf 26} (2009), 185012, 25~pp.;
  {\tt ar$\chi$iv:0906.5477}.

\bibitem{poly1}
B.~Morariu and A.~P.~Polychronakos,
  {\em Finite noncommutative Chern--Simons with a Wilson line and the quantum Hall effect},
  J. High Energy Phys. {\bf 0107} (2001), Paper~6, 15~pp.;
  {\tt ar$\chi$iv:hep-th/0106072}.

\bibitem{fln2}
K.~Noui,
{\em A Model for the motion of a particle in a quantum background},
  Phys.\ Rev.~D {\bf 78} (2008), 105008, 11~pp.;
  {\tt ar$\chi$iv:0807.0969}.


\bibitem{gft2}
D.~Oriti,
 {\em Quantum gravity as a group field theory: a sketch},
  J.\ Phys.\ Conf.\ Ser.\  {\bf 33} (2006), 271--278.

\bibitem{gft3}
 D.~Oriti,
{\em The group field theory approach to quantum gravity}, 
expanded version of an article of the same title in: 
``Approaches to Quantum Gravity --- 
toward a new understanding of space, time, and matter", 
D. Oriti (ed.), Cambridge University Press, Cambridge, 2009,
pp.~310--331;
 {\tt ar$\chi$iv:gr-qc/0607032}.

\bibitem{gft4}
D.~Oriti,
{\em The group field theory approach to quantum gravity: some recent
  results}, extended version of
``Recent progress in group field theory," in: The Planck Scale,
Proceedings of the XXV-th Max Born Symposium, Wroc\l av, 2009,
 J.~Kowalski-Glikman, R.~Durka, M.~Szczachor (eds.),
AIP Conf. Proc., vol.~1196, 2009, pp.~209--218;
 {\tt ar$\chi$iv:0912.2441}.

\bibitem{peskin}
M. E. Peskin and D. V. Schroeder,
{\em An Introduction to Quantum Field Theory}, 
Addison-Wesley, Reading, Mass., 1995.

\bibitem{poly2}
A.~P.~Polychronakos,
{\em Quantum Hall states as matrix Chern--Simons theory},
  J. High Energy Phys. {\bf 0104}  (2001), Paper~11, 20~pp.;
  {\tt ar$\chi$iv:hep-th/0103013}.

\bibitem{book-riv}
V. Rivasseau, {\em From perturbative to constructive renormalization}, 
Princeton University Press, Princeton, 1992.

\bibitem{review-riv} 
  V.~Rivasseau,
{\em An introduction to renormalization},
  in: Poincar\'e Seminar 2002,
B.~Duplantier et al. (eds.), 
Prog. Math. Phys., vol.~30, Birkh\"auser, Basel, 2003, 
pp.~139--177.


\bibitem{rev-riv}
  V.~Rivasseau,
{\em Non-commutative renormalization}, 
in: Quantum spaces, Prog. Math. Phys., vol.~53, 
Birkh\"auser, Basel, 2007, pp.~19--107;
  {\tt ar$\chi$iv:0705.0705}.

\bibitem{csw2}
N.~Seiberg and E.~Witten,
{\em String theory and noncommutative geometry},
  J. High Energy Phys. {\bf 9909} (1999), Paper~32, 93~pp.;
  {\tt ar$\chi$iv:hep-th/9908142}.

\bibitem{sokal}
A. Sokal, {\em The multivariate Tutte polynomial (alias Potts model) for 
graphs and matroids}, in: Surveys in Combinatorics, 2005,
Bridget S.~Webb (ed.), Cambridge University Press,
Cambridge, 2005,
pp.~173--226. 


\bibitem{gauge2}
W. van Suijlekom, Commun. Math. Phys. {\bf 276} (2007), 773--798.

\bibitem{ultimele8}
  A.~Tanasa,
 {\em Algebraic structures in quantum gravity},
  Class.\ Quant.\ Grav.\  {\bf 27} (2010), 095008, 18~pp.;
  {\tt ar$\chi$iv:0909.5631}.


\bibitem{io-param}
A. Tanasa, {\em Parametric representation of a translation-invariant
  renormalizable noncommutative model}, J. Phys.~A {\bf 42} (2010), 
365208, 18~pp.; 
{\tt ar$\chi$iv:0807.2779}.

\bibitem{gp2}
A.~Tanasa,
  {\em Generalization of the Bollob\'as--Riordan polynomial for tensor
    graphs},  J. Math. Phys. {\bf 52} (2011), 073514, 17~pp.;
    {\tt ar$\chi$iv:1012.1798}.


\bibitem{io-kreimer}
A. Tanasa and D. Kreimer, {\em Combinatorial Dyson--Schwinger equation
  in noncommutative field theory},  J. Noncomm. Geom. (in press);
{\tt ar$\chi$iv:0907.2182}.



\bibitem{fab}
 A.~Tanasa and F.~Vignes-Tourneret,
 {\em Hopf algebra of non-commutative field theory},
J.\ Noncomm.\ Geom. {\bf 2} (2008), 125--139;
{\tt ar$\chi$iv:0707.4143}.


\bibitem{tutte1}
W. T. Tutte, {\em Graph Theory}, Addison-Wesley, Reading, Mass.,  1984. 



\bibitem{fab-br}
F. Vignes-Tourneret, {\em The multivariate signed Bollob\'as--Riordan
  polynomial}, Discrete Math. {\bf 309} (2009), 5968--5981;
{\tt ar$\chi$iv:0811.1584}.

\bibitem{mw}
E. W. Weisstein, "Pseudograph." From MathWorld -- A Wolfram Web Resource;
available at
{\tt http://mathworld.wolfram.com/Pseudograph.html}.

\bibitem{teza-yeats}
K.~A.~Yeats,
  {\em Growth estimates for Dyson--Schwinger equations}, Ph.D. thesis, 
Boston Univ., 2008;
{\tt  ar$\chi$iv:0810.2249}.
























\end{thebibliography}
\end{document}